\documentclass[11.5pt]{amsart}
\usepackage{pdfsync}
\usepackage{epsfig}
\usepackage{amsmath, amssymb, euscript,  enumerate}
\usepackage{pxfonts}
\usepackage{bbm}
\usepackage{color}
\usepackage{verbatim}

\newcommand{\supp}{\text{\rm supp}}

\newcommand{\ap}{\alpha}             
\newcommand{\bt}{\beta}
\newcommand{\gm}{\gamma}             \newcommand{\Gm}{\Gamma}
\newcommand{\dt}{\delta}             \newcommand{\Dt}{\Delta}
\newcommand{\vep}{\varepsilon}
\newcommand{\zt}{\zeta}

\newcommand{\ld}{\lambda}            \newcommand{\Ld}{\Lambda}
             
\newcommand{\vp}{\varphi}
\newcommand{\om}{\omega}             \newcommand{\Om}{\Omega}
            \newcommand{\iy}{\infty}
            
\newcommand{\f}{\frac}

\newcommand{\fM}{{\mathfrak M}}

\newcommand{\fR}{{\mathfrak R}}

\newcommand{\BN}{{\mathbb N}}

\newcommand{\BP}{{\mathbb P}}

\newcommand{\BR}{{\mathbb R}}

\newcommand{\cA}{{\mathcal A}}

\newcommand{\cF}{{\mathcal F}}

\newcommand{\cK}{{\mathcal K}}

\newcommand{\cM}{{\mathcal M}}

\newcommand{\cP}{{\mathcal P}}

\newcommand{\cS}{{\mathcal S}}
\newcommand{\cT}{{\mathcal T}}

\newcommand{\cZ}{{\mathcal Z}}

\newcommand{\la}{\langle}          \newcommand{\ra}{\rangle}
\newcommand{\s}{\setminus}         
\newcommand{\n}{\nabla}            \newcommand{\e}{\eta}
\newcommand{\pa}{\partial}        
       
    \newcommand{\ds}{\displaystyle}
\newcommand{\ls}{\lesssim } \newcommand{\gs}{\gtrsim}

 \newcommand{\pf }{\noindent{\it Proof. }}
\newcommand{\rk }{\noindent{\bf Remark. }}
 \newcommand{\bv }{\bar v}

\newcommand{\aee }{\text{\rm a.e.}} 
  \newcommand{\pv }{\text{\rm p.v.}}
\newcommand{\osc }{\text{\rm osc}}

\newcommand{\rL }{{\text{\rm L}}}

\newcommand{\rX}{{\text{\rm X}}}  \newcommand{\rY}{{\text{\rm Y}}}

\newcommand{\bmo}{{\text{\rm BMO}}}
\newcommand{\loc}{{\text{\rm loc}}}
\newcommand{\rh }{{\text{\rm RH}}}
\newcommand{\Div}{{\text{\rm div}}}

\newtheorem{thm}[subsection]{Theorem}
\newtheorem{lemma}[subsection]{Lemma}
\newtheorem{cor}[subsection]{Corollary}

\newtheorem{prop}[subsection]{Proposition}

\newtheorem{defn}[subsection]{Definition}

\numberwithin{equation}{section}

\title[ H\"older continuity of weak solutions ]{ H\"older regularity of weak solutions \\ to nonlocal $p$-Laplacian type Schr\"odinger equations with $A_1^p$-Muckenhoupt potentials }
\author{ Yong-Cheol Kim }
\begin{document}
\begin{abstract} In this paper, using the De Giorgi-Nash-Moser method, we obtain an interior H\"older continuity of weak solutions to nonlocal $p$-Laplacian type Schr\"odinger equations given
by an integro-differential operator $\rL^p_K$ ($\,p>1$) as follows;
\begin{equation*}\begin{cases} \rL^p_K u+V|u|^{p-2} u=0 &\text{ in $\Om$, }\\
                u=g &\text{ in $\BR^n\s\Om$ } \\
                              \end{cases}\end{equation*}
where $V=V_+-V_-$ with $(V_-,V_+)\in L^1_{\loc}(\BR^n)\times L^q_{\loc}(\BR^n)$ for $q>\f{n}{ps}>1$ and $0<s<1$ is a potential such that $(V_-,V_+^{b,i})$ belongs to the $(A_1,A_1)$-Muckenhoupt class and $V_+^{b,i}$ is in the $A_1$-Muckenhoupt class for all $i\in\BN$ ( here, $V_+^{b,i}:=V_+\max\{b,1/i\}/b$ for an almost everywhere positive bounded function $b$ on $\BR^n$ with $V_+/b\in L^q_{\loc}(\BR^n)$, $g\in W^{s,p}(\BR^n)$ and $\Om\subset\BR^n$ is a bounded domain with Lipschitz boundary.) In addition, we get the local boundedness of weak subsolutions of the nonlocal $p$-Laplacian type Schr\"odinger equations. In a different way from \cite{DKP1}, we obtain the logarithmic estimate of the weak supersolutions which play a crucial role in proving the H\"older regularity of the weak solutions.

In particular, we note that all the above results are still working for any nonnegative potential in $L^q_{\loc}(\BR^n)$ ($\,q>\f{n}{ps}>1,  0<s<1$).
\end{abstract}

\address{$\bullet$ Yong-Cheol Kim : Department of Mathematics Education, Korea University, Seoul 02841, Republic of Korea  $\,\&\,$ School of Mathematics, Korea Institute for Advanced Study, Seoul 02455, Republic of Korea }

\email{ychkim@korea.ac.kr}

\maketitle

\tableofcontents

\section{Introduction}
The recent research on fractional and nonlocal partial differential equations has extensively and actively been performed not only in pure mathematics but also in scientific areas that require its concrete applications.
This type of problems appear in various applications such as continuum mechanics, phase transition phenomena that are related to a nonlocal version of the classical Allen-Cahn equation, population dynamics, nonlocal minimal surfaces, a nonlocal version of Schr\"odinger equations for standing waves ( see \cite{BV,CK1,CK2} ), game theory and also constrained variational problems with fractional diffusion arising in the quasi-geostrophic flow model, anomalous diffusions and American options with jump processes ( see \cite{CV,C,S} ).

\,For $q>\f{n}{ps}>1$ ( $p>1$, $0<s<1$), let $\cP^{s,p}_q(\BR^n)$ be the class of all potentials $V=V_+-V_-$ such that (i) $V_-\in L^1_{\loc}(\BR^n)$, (ii) $V_+\in L^q_{\loc}(\BR^n)$, (iii) there is an almost everywhere positive bounded function $b$ on $\BR^n$ so that 

$V_+/b\in L^q_{\loc}(\BR^n)$, $(V_-,V_+^{b,i})$ belongs to the $(A_1,A_1)$-Muckenhoupt class and $V_+^{b,i}$ is in the $A_1$-Muckenhoupt class for all $i\in\BN$, where
$$V^{b,i}_+:=\f{\max\{b,1/i\}}{b}\,V_+.$$
If $V\in\cP^{s,p}_q(\BR^n)$ for $q>\f{n}{ps}>1$ ( $p>1$, $0<s<1$), then we say that 
$V$ is a {\it $A_1^p$-Muckenhoupt potential.} When $p=2$, we call it {\it $A_1$-Muckenhoupt potential.}

\,The aim of this paper is to establish an interior H\"older regularity of weak solutions of nonlocal $p$-Laplacian type Schr\"odinger equations with $A_1^p$-Muckenhoupt potentials and to additionally obtain the local boundedness of weak subsolutions of the nonlocal equation. 

For $p>1$, let $\cK_p$ be the collection of all positive symmetric kernels
satisfying the uniformly ellipticity assumption
\begin{equation}\frac{c_{n,p,s}\,\lambda}{|y|^{n+ps}}\leq K(y)=K(-y)\leq
\frac{c_{n,p,s}\,\Lambda}{|y|^{n+ps}},\,\,0<s<1,\,\,y\in\BR^n\s\{0\},
\end{equation} where $c_{n,p,s}>0$ is the normalization constant given by
\begin{equation*}c_{n,p,s}=\f{\Gm(\f{n+p}{2})\,p\,(1-s)}{\pi^{\f{n-1}{2}}\,\Gm(\f{p+1}{2})}.
\end{equation*}
For $K\in\cK_p$ ($\,p>1$), we consider integro-differential operators $\rL^p_K$ given
by
\begin{equation}\rL^p_K u(x,t)=\pv\int_{\BR^n}H_p(u(x)-u(y))K(x-y)\,dy
\end{equation} where $H_p(t)=|t|^{p-2} t$ for $t\in\BR$. If $p=2$, then we write $\rL^p_K=\rL_K$.
In particular, if $K(y)=c_{n,p,s}|y|^{-n-ps}$,
then $\rL^p_K=(-\Delta)^s_p$ is the fractional $p$-Laplacian and it
is well-known \cite{IN} that $$\lim_{s\to 1^-}(-\Delta)_p^s u=-\Delta_p u$$
for any function $u$ in the Schwartz space $\cS(\BR^n)$, where 
$$\Dt_p u=\Div(|\n u|^{p-2}\n u)$$ is the classical $p$-Laplacian.

\,We are interested in the Dirichlet problem for the nonlocal $p$-Laplacian type Schr\"odinger equation
\begin{equation}\begin{cases} \rL^p_K u+V |u|^{p-2}u=0 &\text{ in $\Om$ }\\
                u=g &\text{ in $\BR^n\s\Om$ } 
                                \end{cases}\end{equation}
where $V\in\cP^{s,p}_q(\BR^n)$ for $q>\f{n}{ps}>1$ ($\,p>1$ and $0<s<1$),
$g\in W^{s,p}(\BR^n)$ and $\Om\subset\BR^n$ is a bounded domain with Lipschitz
boundary. The existence and uniqueness of weak solution to the above nonlocal equation was obtained in \cite{K3} by applying standard technique of calculus of variations.
More precisely speaking about the problem, by employing the De Giorgi-Nash-Moser theory, we obtain an interior H\"older regularity of weak solutions to the nonlocal $p$-Laplacian type Schr\"odinger equations, and also we get the local boundedness of weak subsolutions of the nonlocal equation. Here, we note that the boundary condition is imposed on $\BR^n\s\Om$ with nonlocality. In fact, from the probabilistic point of view, it conforms to the natural phenomenon that a discontinuous L\'evy process on the domain $\Om$ can exit $\Om$ for the first time jumping to any point in $\BR^n\s\Om$.

When $p=2$, the research on the above nonlocal equations was strongly motivated by the study of standing wave solutions of the form
$$\Psi(x,t)=e^{-i\om t}u(x)$$
of the time-dependent nonlocal Schr\"odinger equations
$$i\,\f{\pa\Psi}{\pa t}=\rL_K\Psi+V(x)\Psi$$
which is a fundamental equation of fractional quantum mechanics and fractional quantum physics. This equation was used for the first time in the literature by Laskin ( see \cite{L} ).

It turns out in Section 3 that any potential $V$ in $\cP^{s,p}_q(\BR^n)$ satisfies
\begin{equation}\int_{\BR^n}|\vp(y)|\,V_-(y)\,dy\le\int_{\BR^n}|\vp(y)|\,V_+(y)\,dy
\end{equation}
and 
\begin{equation}\int_{\BR^n}|\vp(y)|^p\,V_-(y)\,dy\le\int_{\BR^n}|\vp(y)|^p\,V_+(y)\,dy
\end{equation}
for every $\vp\in\rY^{s,p}_0(\Om)$, whenever $q>\f{n}{ps}>1$ ($\,p>1$ and $0<s<1$).
As a matter of fact, the inequality (1.4) is a useful tool for the proof of nonlocal Caccioppoli type inequality to be given in Theorem 1.5, and also the inequality (1.5) makes it possible to prove in Lemma 3.3 below that $\rY^{s,p}_0(\Om)$
is a quasi-Banach space.

\,

\,{\bf Notations.} We introduce the notations briefly for the readers as follows.

\,\noindent$\bullet$ For $r>0$, $x_0\in\BR^n$ and $s\in(0,1)$, let us denote by $B_r^0=B_r(x_0)$, $B_r=B_r(0)$. For $n\in\BN$, denote by $|S^{n-1}|$ the $(n-1)$-dimensional surface measure on the unit sphere $S^{n-1}$ of $\BR^n$.

\,\noindent$\bullet$ For two quantities $a$ and $b$, we write $a\lesssim b$ (resp.
$a\gtrsim b$) if there is a universal constant $C>0$ ({\it depending
only on $\ld,\Ld,n,p,s$ and $\Om$}$)$  such that $a\le C\,b$ (resp. $b\le
C\,a$).

\,\noindent$\bullet$ For $a,b\in\BR$, we denote by 
$$a\vee b=\max\{a,b\}\,\,\text{ and }\,\,a\wedge b=\min\{a,b\}.$$

\,\noindent$\bullet$ Let $\cF^n$ be the family of all real-valued
Lebesgue measurable functions on $\BR^n$. 

\,\noindent$\bullet$ For $u\in C(B^0_r)$, we consider the norm 
$$\|u\|_{C(B^0_r)}=\sup_{x\in B^0_r}|u(x)|.$$ For $\gm\in(0,1)$, the {\it $\gm^{th}$ H{$\ddot o$}lder seminorm} of $u$ on $B^0_r$ is defined by
$$[u]_{C^{\gm}(B^0_r)}=\sup_{x,y\in B^0_r,\,x\neq y}\f{|u(x)-u(y)|}{|x-y|^{\gm}}
$$ and the {\it $\gm^{th}$ H{$\ddot o$}lder norm} of $u$ on $B^0_r$ is defined by 
$$\|u\|_{C^{\gm}(B^0_r)}=\|u\|_{C(B^0_r)}+[u]_{C^{\gm}(B^0_r)}.$$

\,\noindent$\bullet$ For $x_0\in\Om$, $p>1$ and $r>0$ with $B^0_r\subset\Om$, the {\it nonlocal tails} of the function $u$ in $B^0_r\subset\Om$ is defined by
\begin{equation}\begin{split}
\cT_r(u;x_0)&=\biggl(\f{ps}{|S^{n-1}|}\,r^{ps}\int_{\BR^n\s B_r(x_0)}\f{|u(y)|^{p-1}}{|y-x_0|^{n+ps}}\,dy\biggr)^{\f{1}{p-1}}.
\end{split}\end{equation}

We now state one of our main results which is called the local boundedness of weak subsolutions to the nonlocal $p$-Laplacian type Schr\"odinger equation (1.3), as follows. 

\begin{thm} Let $V\in\cP^{s,p}_q(\BR^n)$, $g\in W^{s,p}(\BR^n)$ for $q>\f{n}{ps}>1$ $($\,$p>1$, $s\in(0,1)$\,$)$ and $B^0_{2r}\subset\Om$. If $u\in\rY^{s,p}_g(\Om)^-$ is a weak subsolution of nonlocal $p$-Laplacian type Schr\"odinger equation $(1.3)$, then there is a constant $C_0>0$ depending only on $n,s,p,\ld,\Ld$ and $\Om$ such that
\begin{equation*}
\sup_{B^0_r}u\le\dt\,\cT_r(u_+;x_0)+C_0\,\dt^{-\f{(p-1)n}{sp^2}}\biggl(\,\fint_{B^0_{2r}}u_+^p\,dx\biggr)^{\f{1}{p}}
\end{equation*} for any $\dt\in(0,1]$.
\end{thm}

\noindent{\bf{Remark.}} (a) If $u\in\rY^{s,p}_g(\Om)^-$ is a weak subsolution of the nonlocal $p$-Laplacian type Schr${\ddot {\rm o}}$dinger equation (1.4) and $g\in W^{s,p}(\BR^n)$ for $s\in(0,1)$, then we see that $u\in L^p(\Om)$ and $u\le g$ on $\BR^n\s\Om$, and thus $u_+\le g_+$ there. Then it follows from H\"older's inequality and fractional Sobolev inequalities (2.4) and (2.5) that
\begin{equation*}\begin{split}
[\cT_r(u_+;x_0)]^{p-1}&\le\f{ps\,r^{ps}}{|S^{n-1}|}\biggl(\,\int_{\BR^n\s\Om}\f{g_+^{p-1}(y)}{|y-x_0|^{n+ps}}\,dy+\int_{\Om\s B_r(x_0)}\f{|u(y)|^{p-1}}{|y-x_0|^{n+ps}}\,dy\biggr) \\
&\ls\f{ps\,r^{-n(1-\f{1}{p})}}{(((p-1)n+p^2 s)|S^{n-1}|)^{1/p}}\bigl(\,\|g\|^{p-1}_{W^{s,p}(\BR^n)}+\|u\|^{p-1}_{L^p(\Om)}\bigr)<\iy.
\end{split}\end{equation*}

(b) If $u\in\rY^{s,p}_g(\Om)^+$ is a weak supersolution of the nonlocal $p$-Laplacian type Schr${\ddot {\rm o}}$dinger equation (1.4) and $g\in W^{s,p}(\BR^n)$ for $s\in(0,1)$, then $-u$ is its weak subsolution and $u\ge g$ on $\BR^n\s\Om$, and so $u_-\le g_-$ there. Then, as in the above (a), we obtain that
\begin{equation*}[\cT_r(u_-;x_0)]^{p-1}
\ls\f{ps\,r^{-n(1-\f{1}{p})}}{(((p-1)n+p^2 s)|S^{n-1}|)^{1/p}}\bigl(\,\|g\|^{p-1}_{W^{s,p}(\BR^n)}+\|u\|^{p-1}_{L^p(\Om)}\bigr)<\iy.
\end{equation*}
Then, from Theorem 1.1, we easily have that
\begin{equation*}
-\inf_{B^0_r}u=\sup_{B^0_r}(-u)\le\dt\,\cT_r(u_-;x_0)+C_0\,\dt^{-\f{(p-1)n}{sp^2}}\biggl(\,\fint_{B^0_{2r}}u_-^p\,dx\biggr)^{\f{1}{p}}
\end{equation*} for any $\dt\in(0,1]$.

(c)  If $u\in\rY^{s,p}_g(\Om)^+$ is a weak solution of the nonlocal $p$-Laplacian type Schr${\ddot {\rm o}}$dinger equation (1.4) and $g\in W^{s,p}(\BR^n)$ for $s\in(0,1)$, then it follows from (a), (b) and Theorem 1.1 that
\begin{equation*}
\underset{B^0_r}\osc \,u\le 2\dt\,\cT_r(u;x_0)+2C_0\,\dt^{-\f{(p-1)n}{sp^2}}\biggl(\,\fint_{B^0_{2r}}|u|^p\,dx\biggr)^{\f{1}{p}}
\end{equation*} for any $\dt\in(0,1]$.

\,

\,\,\,The following {\it logarithmic estimate} plays a crucial role in proving the H\"older regularity of weak solutions to the nonlocal  $p$-Laplacian type Schr\"odinger equation and in showing that the logarithm of such weak solution is a function with locally bounded mean oscillation. In a different way from \cite{DKP1}, we obtain the logarithmic estimate. We now state it as follows.



\begin{thm} Let $V\in\cP^{s,p}_q(\BR^n)$ and $g\in W^{s,p}(\BR^n)$ for $q>\f{n}{ps}>1$ $($$\,p>1$ and $0<s<1$$)$.
If $u\in\rY^{s,p}_g(\Om)^+$ is a weak supersolution of nonlocal $p$-Laplacian type Schr\"odinger equation $(1.3)$ with $u\ge 0$ in $B^0_R\subset\Om$, then there is a constant $c_0>0$ depending only on $n,s,p,\ld,\Ld$ and $\Om$ such that
\begin{equation*}\begin{split}&\iint_{B^0_r\times B^0_r}\,\biggl|\ln\biggl(\f{u(x)+b}{u(y)+b}\biggr)\biggr|^p\,d_K(x,y)  \\
&\qquad\qquad\qquad\qquad\le c_0\,r^{n-ps}\biggl[\f{1}{b^{p-1}}\biggl(\f{r}{R}\biggr)^{ps}[\cT_R(u_-;x_0)]^{p-1}+\bigl(1+\|V_+\|_{L^q(\Om)}\bigr)\biggr]
\end{split}\end{equation*}for any $b\in(0,1)$ and $r\in(0,R/2)$, where $d_K(x,y)=K(x-y)\,dx\,dy$.
\end{thm}

Employing the De Giorgi-Nash-Moser theory and using Theorem 1.1 and 1.2, we obtain the following H\"older continuity of weak solutions to the nonlocal  $p$-Laplacian type Schr\"odinger equation, and also we can easily derive Corollary 1.4 as a natural by-product of Theorem 1.3.

\begin{thm} Let $V\in\cP^{s,p}_q(\BR^n)$, $g\in W^{s,p}(\BR^n)$ for $q>\f{n}{ps}>1$ $($$\,p>1$, $0<s<1$$)$, and let $B^0_{2R}\subset\Om$.
If $u\in\rY^{s,p}_g(\Om)$ is a weak solution of the nonlocal $p$-Laplacian type Schr\"odinger equation $(1.3)$, then there exist constants $\e_0^-\in(0,\f{ps}{2(p-1)})$ and $\e_0^+\in(\f{ps}{2(p-1)},\f{ps}{p-1})$ such that $u$ is locally $\e$-H\"older continuous in $\Om$ for any $\e\in(0,\e_0^-]\cup[\e_0^+,\f{ps}{p-1})$. Furthermore, for each $x_0\in\Om$ and for each $\e\in(0,\e_0^-]\cup\bigl[\e_0^+,\f{ps}{p-1}\bigr)$, we have that
\begin{equation}\underset{B^0_r}\osc\,u\ls\biggl(\f{r}{R}\biggr)^{\e}\,\biggl[\cT_R(u;x_0)+\biggl(\,\fint_{B^0_{2R}}|u(x)|^p\,dx\biggr)^{\f{1}{p}}\biggr]
\end{equation}
for any $r\in (0,R/2)$. Here it turns out that there exist universal constants $c_0, c_*>0$ such that
$$\e_0^{\pm}=\f{\ln\biggl(\,\ds\f{1\pm\sqrt{1-4\,\dt^{\f{ps}{p-1}}}}{2}\,\biggr)}{\ln\dt}\,\,\text{ for $\dt=e^{-(c_0/c_*)(1+\|V_+\|_{L^q(\Om)})^{1/p}}\wedge\bigl(\f{1}{4}\bigr)^{\f{p-1}{ps}}$. }$$
\end{thm}

\,

\,\,\,The next corollary can easily be obtained by applying Theorem 1.3 and employing the interpolation on H\"older spaces between $C^{\e_0^-}(B^0_r)$ and $C^{\e_0^+}(B^0_r)$, which eventually fill up an interior $\e$-H\"older continuity of $u$ in $\Om$ for all $\e\in(\e_0^-,\e_0^+)$.

\begin{cor} Let $V\in\cP^{s,p}_q(\BR^n)$, $g\in W^{s,p}(\BR^n)$ for $q>\f{n}{ps}>1$ $($$\,p>1$, $0<s<1$$)$, and let $B^0_{2R}\subset\Om$.
If $u\in\rY^{s,p}_g(\Om)$ is a weak solution of the nonlocal $p$-Laplacian type Schr\"odinger equation $(1.3)$, then we have the following estimate
\begin{equation}
\sup_{r\in(0,R/2)}\,\|u\|_{C^{\e}(B^0_r)}\ls \f{1}{R^\e}\,\biggl[\cT_{R}(u;x_0)+\biggl(\,\fint_{B^0_{2R}}|u(x)|^p\,dx\biggr)^{\f{1}{p}}\biggr]
\end{equation} for any $\e\in \bigl(0,\f{ps}{p-1}\bigr)$.
\end{cor}

\rk If $\f{p-1}{p}<s<1$, then we can expect the better regularity, i.e. $C^{1,\ap}$-estimate for some $\ap\in(0,1)$.

\,

\,\,\,As a basic tool for our main results, we show that any weak subsolution of the nonlocal  $p$-Laplacian type Schr\"odinger equations enjoys the following {\it nonlocal Caccioppoli type inequality.}

\begin{thm} Let $V\in\cP^{s,p}_q(\BR^n)$, $g\in W^{s,p}(\BR^n)$ for $q>\f{n}{ps}>1$ $($$\,p>1$, $s\in(0,1)$$)$, and let $B^0_{2r}\subset\Om$.
If $u\in\rY^{s,p}_g(\Om)^-$ is a weak subsolution of the nonlocal $p$-Lapacian type Schr\"odinger  equation $(1.3)$, then for any nonnegative $\zt\in C_c^{\iy}(B^0_r)$ we have the following estimate
\begin{equation*}\begin{split}
&\int_{B^0_r}[w(y)\zt(y)]^pV(y)\,dy+\iint_{B^0_r\times B^0_r}|\zt(x)w(x)-\zt(y)w(y)|^p\,d_K(x,y) \\
&\qquad\qquad\le 2^{2p+1}\bigl(\f{1}{4}+c_p\bigr)\iint_{B^0_r\times B^0_r}[w(x)\vee w(y)]^p|\zt(x)-\zt(y)|^p\,d_K(x,y) \\
&\qquad\qquad\qquad+2^{p+2}\biggl(\,\sup_{x\in\supp(\zt)}\int_{\BR^n\s B^0_r}w^{p-1}(y)\,K(x-y)\,dy\biggr)\|w\zt^p\|_{L^1(B^0_r)}
\end{split}\end{equation*} where $w=(u-M)_+$ 
for $M\in(0,\iy)$ and $c_p=\f{1}{2}[2(p-1)]^{p-1}$.
\end{thm}

\rk (a) When $p=2$ and $s=1$, the study of the classical Schr\"odinger operator, i.e. local Schr\"odinger operator $-\Dt+V$ has been ongoing actively and widely in analysis area in Mathematics and Mathematical Physics (refer to \cite{AS, CFG, F, S1, S2, Si}).

\,(b) In case that $p=2$ and $V\in L^q_{\loc}(\BR^n)$ with $q>\f{n}{2s}$ ($0<s<1$) is nonnegative, it is known in \cite{CK1} that a fundamental solution for nonlocal Schr\"odinger operator $\rL_K+V$ exists and its decay can be obtained. Also, under an additional restiction that the potential $V$ is in a reverse H\"older class $\rh^\tau$ for $\tau>\f{n}{2s}>1$ ($0<s<1$), the $L^a-L^b$ estimate for the Schr\"odinger operator $\rL_K+V$ was obtained inside certain trapezoidal region $\cZ$ which is consist of $(\f{1}{a},\f{1}{b})$ and also the weak type $L^a-L^b$ estimate was partially obtained on the boundary of the region $\cZ$ (see \cite{CK2}).

\,(c) When $p=2$, $0<s<1$ and $V$ is an $A_1$-Muckenhoupt potential, it was shown in \cite{K3} that a fundamental solution for nonlocal Schr\"odinger operator $\rL_K+V$ exists and its decay can be obtained. Also, H\"older continuity and nonlocal Harnack inequalities for $\rL_K+V$ were obtained in \cite{K1} and \cite{K2}.

\,(d) When $V=0$ and $0<s<1$, the result of this problem was obtained by Di Castro, Kuusi and Palatucci \cite{DKP}; as a matter of fact, when $p\in(1,\iy)$, they proved nonlocal Harnack inequalities for elliptic nonlocal $p$-Laplacian equations there, and also they obtained H$\ddot {\rm o}$lder regularity in \cite{DKP1}.

\,(e) When $p=2$ and $0<s<1$, nonlocal Harnack inequalities for locally nonnegative weak solutions of nonlocal heat equations was obtained in \cite{K} by applying the De Giorgi-Nash-Moser theory.

\,\, The paper is organized as follows. In Section 2, we furnish the
function spaces and the definition of weak solutions of the nonlocal
Schr\"odinger equation given in (1.3), and also give a well-known
lemma which is useful in applying the De Giorgi-Nash-Moser theory. 
In Section 3, we give a brief introduction about weighted norm inequalities and the $A_p$-Muckenhoupt class. Additionally, we furnish several examples about sign-changing potentials in the class $\cP^{s,p}_q(\BR^n)$ ($q>\f{n}{ps}>1$, $p>1$, $0<s<1$).
In Section 4, we obtain a sort of {\it nonlocal Caccioppoli type inequality} and several useful local properties of weak solutions to the nonlocal $p$-Laplacian type Schr\"odinger equation by using it. In Section 5, we show that the logarithm of a weak solution to the nonlocal $p$-Laplacian type Schr\"odinger equation becomes a function with locally bounded mean oscillation. In Section 6, we get an interior H\"older continuity of weak solutions to the nonlocal  $p$-Laplacian type Schr\"odinger equation by applying the results obtained in Section 4 and Section 5.

\section{Preliminaries} 
Let $\Om\subset\BR^n$ be a bounded domain with Lipschitz boundary
and let $K\in\cK_p$ for $p>1$. For $p>1$ and $0<s<1$, let $\rX^{s,p}(\Om)$ be the linear function space of all
Lebesgue measurable functions $v\in\cF^n$ such that $v|_\Om\in L^p(\Om)$
and
\begin{equation*}\iint_{\BR^{2n}_\Om}\f{|v(x)-v(y)|^p}{|x-y|^{n+ps}}\,dx\,dy<\iy
\end{equation*}
where $\BR^{2n}_S:=\BR^{2n}\s(S^c\times S^c)$ for a set $S\subset\BR^n$. We
also set
\begin{equation}\rX^{s,p}_0(\Om)=\{v\in\rX^{s,p}(\Om):v=0\,\,\aee\text{ in $\BR^n\s\Om$ }\}
\end{equation}
Since $C^2_0(\Om)\subset\rX^{s,p}_0(\Om)$, we see that $\rX^{s,p}(\Om)$ and $\rX^{s,p}_0(\Om)$ are nonempty. Then we see that $(\rX^{s,p}(\Om),\|\cdot\|_{\rX^{s,p}(\Om)})$ is a normed space with the norm $\|\cdot\|_{\rX^{s,p}(\Om)}$ given by
\begin{equation}\|v\|_{\rX^{s,p}(\Om)}=\|v\|_{L^p(\Om)}+\biggl(\iint_{\BR^{2n}_\Om}\f{|v(x)-v(y)|^p}{|x-y|^{n+ps}}\,dx\,dy\biggr)^{\f{1}{p}}<\iy
\end{equation} for $v\in\rX^{s,p}(\Om)$.
For $p\ge 1$, we denote by $W^{s,p}(\Om)$ the usual fractional Sobolev space
with the norm
\begin{equation}\|v\|_{W^{s,p}(\Om)}:=\|v\|_{L^p(\Om)}+[v]_{W^{s,p}(\Om)}<\iy
\end{equation} where the seminorm $[\,\cdot\,]_{W^{s,p}(\Om)}$ is defined by
$$[v]_{W^{s,p}(\Om)}=\biggl(\iint_{\Om\times\Om}
\f{|v(x)-v(y)|^p}{|x-y|^{n+ps}}\,dx\,dy\biggr)^{\f{1}{p}}.$$
When $\Om=\BR^n$ in (2.3), similarly we define the spaces $W^{s,p}(\BR^n)$ for $p\ge 1$ and $s\in(0,1)$.

If $p\ge 1$ and $s\in(0,1)$ satisfy $ps<n$, then it is well-known \cite{DPV} that there exists a universal constant $c=c(n,p,s,\Om)>0$ such that
\begin{equation}\|f\|_{L^\tau(\Om)}\le c\,\|f\|_{W^{s,p}(\Om)}
\end{equation} for any $f\in W^{s,p}(\Om)$ and $\tau\in[p,p_*]$, where $p_*=\f{pn}{n-ps}$. Furthermore, there is a universal constant $c=c(n,p,s)>0$ such that
\begin{equation}\|f\|_{L^{\tau}(\BR^n)}\le c\,\|f\|_{W^{s,p}(\BR^n)},\forall\tau\in[p,p_*]\,\,\text{ and }\,\,\|f\|_{L^{p_*}(\BR^n)}\le c\,[f]_{W^{s,p}(\BR^n)}
\end{equation} for any $f\in W^{s,p}(\BR^n)$.
Using (2.5), we easily see that there exists a constant $c>1$ depending only on $n,p,s$ and $\Om$ such that
\begin{equation}\|u\|_{\rX^{s,p}_0(\Om)}\le\|u\|_{\rX^{s,p}(\Om)}\le c\,\|u\|_{\rX^{s,p}_0(\Om)}
\end{equation}
for any $u\in\rX^{s,p}_0(\Om)$, where
\begin{equation}\|u\|_{\rX^{s,p}_0(\Om)}:=\biggl(\iint_{\BR^{2n}_\Om}\f{|u(x)-u(y)|^p}{|x-y|^{n+ps}}\,dx\,dy\biggr)^{\f{1}{p}}.
\end{equation} Thus $\|\cdot\|_{\rX^{s,p}_0(\Om)}$ is a norm on $\rX^{s,p}_0(\Om)$ which is equivalent to (2.2). 
By using the change of variables, we can easily derive the following version of the {\it fractional Sobolev inequality} (2.4).

\begin{prop} Let $B_R$ be a ball with radius $R>0$. If $s\in(0,1)$ and $p\in[1,\iy)$ with $sp<n$, then there is a constant $c=c(n,p,s)>0$ such that
\begin{equation*}\|f\|_{L^\tau(B_R)}\le c\,R^{-n(\f{1}{p}-\f{1}{\tau})}\|f\|_{L^p(B_R)}+c\,R^{-n(\f{1}{p}-\f{1}{\tau})+s}[f]_{W^{s,p}(B_R)}
\end{equation*} for any $\tau\in[p,p_*]$.
In particular, if $\tau=p_*:=\f{pn}{n-ps}$, then we have that
\begin{equation}\|f\|_{L^{p_*}(B_R)}\le c\,R^{-s}\|f\|_{L^p(B_R)}+c\,[f]_{W^{s,p}(B_R)}.
\end{equation}
\end{prop}

\,

\,\,\,For $g\in W^{s,p}(\BR^n)$, we consider the convex subsets of $\rX^{s,p}(\Om)$ by \begin{equation*}\begin{split}\rX^{s,p}_g(\Om)^\pm&=\{v\in\rX^{s,p}(\Om):(g-v)_{\pm}\in\rX^{s,p}_0(\Om)\},\\
\rX^{s,p}_g(\Om)&:=\rX_g^{s,p}(\Om)^+\cap\rX_g^{s,p}(\Om)^-=\{v\in\rX^{s,p}(\Om):g-v\in\rX^{s,p}_0(\Om)\}.
\end{split}\end{equation*}
For $g\in W^{s,p}(\BR^n)$ and a potential $V\in\cP^{s,p}_q(\BR^n)$ with $q>\f{n}{ps}>1$ ($p>1$, $0<s<1$), let 
$$\rY^{s,p}(\Om)=\rX^{s,p}(\Om)\cap L^p_V(\Om)\,\,\text{ and }\,\, \rY^{s,p}_g(\Om)=\rX^{s,p}_g(\Om)\cap L^p_V(\Om)$$ where $L^p_V(\Om)$ is the weighted $L^p$ class of all real-valued measurable functions $u$ on $\BR^n$ satisfying
\begin{equation*}\begin{split}
-\iy<\|u\|^p_{L^p_V(\Om)}&:=\int_\Om|u(y)|^p\,V_+(y)\,dy-\int_\Om|u(y)|^p\,V_-(y)\,dy \\
&:=\|u\|^p_{L^p_{V_+}(\Om)}-\|u\|^p_{L^p_{V_-}(\Om)}<\iy.
\end{split}\end{equation*}
That is, we see that $u\in L^p_V(\Om)$ if and only if $u\in L^p_{V_+}(\Om)\cap L^p_{V_-}(\Om)$.
Here, we note that $\|u\|^p_{L^p_V(\Om)}$ need not be nonnegative, and so the class $L^p_V(\Om)$ is not always a normed space. 

Also we consider function spaces $\rY^{s,p}_g(\Om)^+$ and $\rY^{s,p}_g(\Om)^-$ defined by
$$\rY^{s,p}_g(\Om)^\pm=\{u\in\rY^{s,p}(\Om):(g-u)_\pm\in\rY^{s,p}_0(\Om)\}.$$
Then we see that 
$$\rY^{s,p}_g(\Om)=\rY^{s,p}_g(\Om)^+\cap\rY^{s,p}_g(\Om)^-.$$ If $u=g=0$ in $\BR^n\s\Om$, then 
we easily know that $\rY^{s,p}_0(\Om)=\rX^{s,p}_0(\Om)\cap L^p_V(\Om)$ need not be a Banach space. However, if $V\in\cP^{s,p}_q(\BR^n)$ for $q>\f{n}{ps}>1$ ($\,p>1$, $0<s<1$), then it turns out in Lemma 3.2 below that the class $\rY^{s,p}_0(\Om)$ is a quasi-Banach space with the quasinorm  $\|\cdot\|_{\rY^{s,p}_0(\Om)}$ given by
$$\|u\|^p_{\rY^{s,p}_0(\Om)}:=\|u\|^p_{\rX^{s,p}_0(\Om)}+\int_\Om|u(y)|^p V(y)\,dy,\,\,u\in\rY^{s,p}_0(\Om),$$
$\rY^{s,p}_0(\Om)=\rX^{s,p}_0(\Om)$ and they are quasinorm-equivalent.

\,In order to define weak solutions of the nonlocal equation (1.3), we consider a bilinear form $\la\cdot,\cdot\ra_{H_p,K}:\rX^{s,p}(\Om)\times\rX^{s,p}(\Om)\to\BR$ defined by
$$\la u,v\ra_{H_p,K}=\iint_{\BR^n\times\BR^n}H_p(u(x)-u(y))(v(x)-v(y))\,d_K(x,y)$$ where $d_K(x,y):=K(x-y)\,dx\,dy$.

\begin{defn} Let $V\in\cP^{s,p}_q(\BR^n)$ and $g\in W^{s,p}(\BR^n)$ for $q>\f{n}{ps}>1$ $($$\,p>1$, $0<s<1$$)$. Then we say that a function $u\in\rY^{s,p}_g(\Om)^-$ $(\,u\in\rY^{s,p}_g(\Om)^+\,)$ is a {\rm weak subsolution (\,weak supersolution\,)} of the nonlocal $p$-Laplacian type Schr\"odinger equation $(1.3)$, if it satisfies
\begin{equation}\la u,\vp\ra_{H_p,K}+\int_{\BR^n}V(x)|u(x)|^{p-2}u(x)\vp(x)\,dx\le 0\,\,(\,\ge 0\,)\end{equation}
for any nonnegative $\vp\in\rY^{s,p}_0(\Om)$.
Also, we say that a function $u$ is a {\rm weak solution} of the nonlocal equation $(1.3)$, if it is both a weak subsolution and a weak supersolution, i.e.
\begin{equation}\la u,\vp\ra_{H_p,K}+\int_{\BR^n}V(x)|u(x)|^{p-2}u(x)\vp(x)\,dx=0\,\,\,\text{ for any $\vp\in\rY^{s,p}_0(\Om)$. }\end{equation}
\end{defn}

\,\,To prove our results, we need a well-known lemma \cite{GT} that is useful in applying the De Giorgi-Nash-Moser method.

\begin{lemma} Let $\{N_k\}_{k=0}^{\iy}\subset\BR$ be a sequence of positive numbers such that
$$N_{k+1}\le d_0\,e_0^k N_k^{1+\e}\,\,\text{ for all $k\in\BN\cup\{0\}$ }$$
where $d_0,\e>0$ and $e_0>1$. If $N_0\le d_0^{-1/\e} e_0^{-1/\e^2}$, then we have that $N_k\le e_0^{-k/\e}\,N_0$ for any $k=0,1,\cdots$ and moreover $\lim_{k\to\iy}N_k=0$.
\end{lemma}

\,\,We need several elementary inequalities which are useful in proving Theorem 1.2 and Theorem 1.5.

\begin{lemma} $(a)$ If $\ap,\bt\in\BR$ and $A,B\ge 0$, then we have the inequality
\begin{equation*}|\bt-\ap|^{p-2}(\bt-\ap)(\bt B^p-\ap A^p)\ge-c_p\,\bigl(|\ap|+|\bt|\bigr)^p\,|B-A|^p
\end{equation*} for any $p>1$, where $c_p=\f{(p-1)^{p-1}}{2}$.

$(b)$ If $\ap,\bt\in\BR$ with $\bt\ge\ap$ and $A,B\ge 0$, then we have the inequality
\begin{equation*}(\bt-\ap)^{p-1}(\bt B^p-\ap A^p)\ge\f{1}{4}(\bt-\ap)^p(A^p+B^p)-d_p\,\bigl(|\ap|+|\bt|\bigr)^p\,|B-A|^p
\end{equation*} for any $p>1$, where $d_p=\f{1}{2}[2(p-1)]^{p-1}$.

$(c)$ If $A\ge B\ge 0$ and $p>1$, then $(A-B)^{p-1}\ge b_p A^{p-1}-B^{p-1}$ for $p>1$, where 
$b_p=\mathbbm{1}_{(1,2]}(p)+2^{-(p-1)}\mathbbm{1}_{(2,\iy)}(p)$.
\end{lemma}

\pf (a) Without loss of generality, we may assume that $\ap\neq\bt$; otherwise, the inequality is trivial.
By simple calculation, we have that
\begin{equation*}\begin{split}I&:=|\bt-\ap|^{p-2}(\bt-\ap)(\bt B^p-\ap A^p) \\
&=|\bt-\ap|^p B^p+|\bt-\ap|^{p-2}(\bt-\ap)\ap(B^p-A^p) \\
&\ge |\bt-\ap|^p B^p-|\bt-\ap|^{p-1}\,|\ap|\,|B^p-A^p|.
\end{split}\end{equation*}
Similarly, we have the following inequality
\begin{equation*}I\ge |\bt-\ap|^p A^p-|\bt-\ap|^{p-1}\,|\bt|\,|B^p-A^p|.
\end{equation*}
This yields that
\begin{equation}2I\ge |\bt-\ap|^p (A^p+B^p)-|\bt-\ap|^{p-1}\,(|\ap|+|\bt|)\,|B^p-A^p|.
\end{equation}
Without loss of generality assuming that $B\ge A$, then it follows that
\begin{equation}\begin{split}B^p-A^p&=\int_A^B \f{d}{dt}(t^p)\,dt=p\int_A^B t^{p-1}\,dt \\
&\le p\, B^{p-1}\,(B-A)\le p\bigl(A^p+B^p\bigr)^{\f{p-1}{p}}(B-A).
\end{split}\end{equation}
Thus, by (2.12) and Young's inequality with any $\vep>0$, we have that
\begin{equation}\begin{split}|B^p-A^p|&\le p\bigl(A^p+B^p\bigr)^{\f{p-1}{p}}|B-A| \\
&\le p\vep(A^p+B^p)+\f{1}{\bigl(\vep\f{p}{p-1}\bigr)^{p-1}}\,|B-A|^p.
\end{split}\end{equation}
From (2.11) and (2.13), we have that
\begin{equation}\begin{split}2 I&\ge|\bt-\ap|^{p-1}(A^p+B^p)\bigl[|\bt-\ap|-p\vep(|\ap|+|\bt|)\bigr] \\
&\qquad-\f{1}{\bigl(\vep\f{p}{p-1}\bigr)^{p-1}}\,|\bt-\ap|^{p-1}\,(|\ap|+|\bt|)\,|B-A|^p.
\end{split}\end{equation}
Taking $\vep=\ds\f{|\bt-\ap|}{p(|\ap|+|\bt|)}>0$ in (2.14), the required inequality can be obtained.

(b) Without loss of generality, we may assume that $\ap\neq\bt$; for, if not, the inequality is obvious.
By simple calculation, we have that
\begin{equation*}\begin{split}J&:=(\bt-\ap)^{p-1}(\bt B^p-\ap A^p) \\
&=(\bt-\ap)^p B^p+(\bt-\ap)^{p-1}\ap(B^p-A^p) \\
&\ge (\bt-\ap)^p B^p-(\bt-\ap)^{p-1}\,|\ap|\,|B^p-A^p|.
\end{split}\end{equation*}
Similarly, we have the following inequality
\begin{equation*}J\ge (\bt-\ap)^p A^p-(\bt-\ap)^{p-1}\,|\bt|\,|B^p-A^p|.
\end{equation*}
This yields that
\begin{equation}2J\ge (\bt-\ap)^p (A^p+B^p)-(\bt-\ap)^{p-1}\,(|\ap|+|\bt|)\,|B^p-A^p|.
\end{equation}
Assuming $B\ge A$ without loss of generality, then it follows that
\begin{equation}\begin{split}B^p-A^p&=\int_A^B \f{d}{dt}(t^p)\,dt=p\int_A^B t^{p-1}\,dt \\
&\le p\, B^{p-1}\,(B-A)\le p\bigl(A^p+B^p\bigr)^{\f{p-1}{p}}(B-A).
\end{split}\end{equation}
Thus, by (2.16) and Young's inequality with any $\vep>0$, we have that
\begin{equation}\begin{split}|B^p-A^p|&\le p\bigl(A^p+B^p\bigr)^{\f{p-1}{p}}|B-A| \\
&\le p\vep(A^p+B^p)+\f{1}{\bigl(\vep\f{p}{p-1}\bigr)^{p-1}}\,|B-A|^p.
\end{split}\end{equation}
From (2.15) and (2.17), we have that
\begin{equation}\begin{split}2 J&\ge(\bt-\ap)^{p-1}(A^p+B^p)\bigl[(\bt-\ap)-p\vep(|\ap|+|\bt|)\bigr] \\
&\qquad-\f{1}{\bigl(\vep\f{p}{p-1}\bigr)^{p-1}}\,(\bt-\ap)^{p-1}\,(|\ap|+|\bt|)\,|B-A|^p.
\end{split}\end{equation}
Taking $\vep=\ds\f{\bt-\ap}{2p(|\ap|+|\bt|)}>0$ in (2.18), the required inequality can be obtained.

(c) If $a,b\ge 0$, then it easily follows from the fact that 
$(a+b)^{p-1}\le a^{p-1}+b^{p-1}$ for $1<p\le 2$ and $(a+b)^{p-1}\le 2^{p-1}\bigl(a^{p-1}+b^{p-1}\bigr)$ for $p>2$. Set $a=A-B$ and $b=B$.
\qed

\begin{lemma}$($\cite{DKP1}$)$ If $p\ge 1$, $\vep\in(0,1]$ and $a,b\ge 0$, then the inequality holds
$$a^p\le b^p+c_p\vep\,b^p+(1+c_p\vep)\vep^{1-p}|a-b|^p,$$
where $c_p=(p-1)\Gm(1\vee(p-2))$ for the standard Gamma function $\Gm$.
\end{lemma}

\section{Weighted norm inequalities on the $A_p$-Muckenhoupt class }

In this section, we briefly introduce the $A_p$-Muckenhoupt class for $p\ge 1$ and we prove that any potential $V$ in $\cP^{s,p}_q(\BR^n)$ (\,$q>\f{n}{ps}>1$, $p>1$, $0<s<1$) satisfies certain weighted norm inequalities related with $(V_-,V_+)$. 

By a {\it weight} on $\BR^n$ equipped with Lebesgue measure, we mean a locally integrable function on $\BR^n$ taking value in $[0,\iy)$ almost everywhere. For $f\in L^1_{\rm loc}(\BR^n)$ and $x\in\BR^n$, the Hardy-Littlewood maximal function $\cM f$  is defined by
$$\cM f(x)=\sup_{Q^x}\fint_{Q^x}|f(y)|\,dy,$$
where the supremum is taken over all cubes $Q^x$ with center $x$. For a pair of weights $(\upsilon,\om)$, the quantity $[\upsilon,\om]_{(A_p,A_p)}$ is defined by
\begin{equation*}[\upsilon,\om]_{(A_p,A_p)}=\begin{cases}
\ds\sup_{Q}\biggl(\,\fint_{Q} \upsilon(y)\,dy\biggr)\biggl(\,\fint_{Q} \om(y)^{-\f{1}{p-1}}\,dy\biggr)^{p-1}, &1<p<\iy, \\
\ds\sup_{Q}\biggl(\,\fint_{Q} \upsilon(y)\,dy\biggr)\,\|\om^{-1}\|_{L^\iy(Q)}, &p=1,
\end{cases}\end{equation*} where the supremum is taken over all cubes $Q\subset\BR^n$.
For $p\in[1,\iy)$, we say that $(\upsilon,\om)\in (A_p,A_p)$ if $[\upsilon,\om]_{(A_p,A_p)}<\iy$, and $\om\in A_p$ if $[\om]_{A_p}:=[\om,\om]_{(A_p,A_p)}<\iy$. For $1\le p<\iy$, it is well-known in \cite{Gr} that $(A_1,A_1)\subset(A_p,A_p)$, $A_1\subset A_p$ and $(\upsilon,\om)\in (A_p,A_p)$ is equivalent to the mapping property that 
$\cM:L^p_\om(\BR^n)\to L^{p,\iy}_\upsilon(\BR^n)$ is bounded; that is, there is a universal constant $C_{n,p}>0$ such that
\begin{equation}\sup_{t>0}\bigl[t\,\upsilon\bigl(\{x\in\BR^n:\cM f(x)>t\}\bigr)^{\f{1}{p}}\bigr]\le C_{n,p}\biggl(\int_{\BR^n}|f(y)|^p\om(y)\,dy\biggr)^{\f{1}{p}}
\end{equation}
for any $f\in L^p_\om(\BR^n)$. Here, we denote by $\upsilon(E)=\int_E\upsilon(y)\,dy$ for a set $E\subset\BR^n$. The reader can also refer to \cite{Gr} for these concepts in Fourier analysis.

\,We shall now furnish several examples in the class $\cP^{s,p}_q(\BR^n)$ for $q>\f{n}{ps}>1$ ($\,p>1$, $0<s<1$) mentioned in the above introduction:

\,(a) If $h_\e(x)=|x|^\e$ for $\ap\in\BR$, then it is easy to check that 
$$[h_\e]_{A_1}<\iy\,\,\text{ if and only if }\,\,-n<\e\le 0.$$ If we consider a sign-changing potential $$V(x)=|x|^{\e/q}\cos(|x|)$$ with $\e\in(-n,0]$ and $q>\f{n}{ps}>1$ ($p>1$, $0<s<1$), then we can easily check that $V_+\in L^q_{\loc}(\BR^n)$, $V_-\in L^1_{\loc}(\BR^n)$, $V_+/b\in L^q_{\loc}(\BR^n)$, $(V_-,V^{b,i}_+)\in (A_1,A_1)$ and $V^{b,i}_+\in A_1$ for all $i\in\BN$, where $b(x)=\cos^+(|x|)$.

\,(b) Let $\nu$ be a Borel measue on $\BR^n$ satisfying that
$$\fM\nu(x):=\sup_{Q^x}\f{\nu(Q^x)}{|Q^x|}\le C\,\,\text{ for $\aee$ $x\in\BR^n$,}$$
where the supremum is taken over all cubes $Q^x$ with center $x\in\BR^n$.
If we set $g_\gm(x)=[\fM\nu(x)]^{\gm}$ for $\gm\in(0,1)$, then it is known in \cite{GR} that $[g_\gm]_{A_1}<\iy$. We consider the following sign-changing potential 
$$V(x)=g_\gm(x)\sin(1/|x|)$$ for $\gm\in(0,1).$
If $q>\f{n}{ps}>1$ for $p>1$ and $0<s<1$, then it is easy to check that $V_+\in L^q_{\loc}(\BR^n)$, $V_-\in L^1_{\loc}(\BR^n)$, $V_+/b\in L^q_{\loc}(\BR^n)$, $(V_-,V^{b,i}_+)\in(A_1,A_1)$ and $V^{b,i}_+\in A_1$ for all $i\in\BN$, where $b(x)=\sin^+(1/|x|)$.

\,(c) Let $\om(x)=\ln(1/|x|)\chi_{B(0;1/e)}+\chi_{\BR^n\s B(0;1/e)}$. Then it is easy to check that $\om\in A_1$. We consider the following sign-changing potential
$$V(x)=\bigl[\ln(1/|x|)\chi_{B(0;1/e)}+\chi_{\BR^n\s B(0;1/e)}\bigr]\cos(1/|x|).$$ From properties of the Gamma function, if $q>\f{n}{ps}>1$ for $p>1$ and $0<s<1$ then we can easily check that $V_+\in L^q_{\loc}(\BR^n)$, $V_-\in L^1_{\loc}(\BR^n)$, $V_+/b\in L^q_{\loc}(\BR^n)$, $(V_-,V^{b,i}_+)\in (A_1,A_1)$ and $V^{b,i}_+\in A_1$ for all $i\in\BN$, where $b(x)=\cos^+(1/|x|)$. 

\,

\,\,Next, we need several fundamental lemmas which is useful in proving Theorem 3.3.

\begin{lemma} If $V\in\cP^{s,p}_q(\BR^n)$ for $q>\f{n}{ps}>1$, $p>1$ and $0<s<1$, then 
$$\int_{\BR^n}|\vp(y)|^p\,V_-(y)\,dy\le\int_{\BR^n}|\vp(y)|^p\,V_+(y)\,dy$$
for any $\vp\in\rY^{s,p}_0(\Om)$.
\end{lemma}

\pf From the exactly same way as the proof of Theorem 3.6 \cite{K2}, we see that
\begin{equation}\int_{\BR^n}|\vp(y)|^p\,V(y)\,dy\ge 0\end{equation}
for any $\vp\in C_c^\iy(\BR^n)$. 
Take any $\vp\in\rY^{s,p}_0(\Om)$. Since $C^\iy_c(\Om)$ is dense in $\rX^{s,p}_0(\Om)$ (see \cite{FSV} and Theorem 1.4.2.2 in \cite{G}), we can take a sequence $\{\vp_k\}\subset C^\iy_c(\Om)$ such that 
$$\vp_k\to\vp\,\,\text{ in $\rX^{s,p}_0(\Om)$.}$$
Since $V\in\BP^{s,p}_q(\BR^n)$, there is a nonnegative bounded function $b$ on $\BR^n$ such that
\begin{equation}V^b_+:=\f{V_+}{b}\in L^q_{\rm loc}(\BR^n), (V_-,V^{b,i}_+)\in(A_1,A_1)\text{ and }V^{b,i}_+\in A_1
\end{equation} for all $i\in\BN$.
Then we claim that
\begin{equation}\lim_{k\to\iy}\int_{\BR^n}|\vp_k(y)|^p\,V_+(y)\,dy=\int_{\BR^n}|\vp(y)|^p\,V_+(y)\,dy;
\end{equation}
indeed, we note that $\|\vp_k\|_{\rX_0^{s,p}(\Om)}\le 2\|\vp\|_{\rX_0^{s,p}(\Om)}$ for all sufficiently large $k\in\BN$, and also we see that, for any $y\in\BR^n$,
\begin{equation*}\begin{split}
|\vp_k(y)|^p-|\vp(y)|^p&=\int_{|\vp(y)|}^{|\vp_k(y)|}\f{d}{d\tau}\tau^p\,d\tau \\
&\le p\,\bigl(|\vp_k(y)|-|\vp(y)|\bigr)\,\bigl(|\vp_k(y)|\vee|\vp(y)|\bigr)^{p-1}.
\end{split}\end{equation*}
Thus it follows from the fractional Sobolev inequality and H\"older's inequality that
\begin{equation}\begin{split}
&\biggl|\int_{\BR^n}|\vp_k(y)|^p\,V_+(y)\,dy-\int_{\BR^n}|\vp(y)|^p\,V_+(y)\,dy\biggr| \\
&\qquad\qquad\le\int_{\BR^n}\bigl||\vp_k(y)|^p-|\vp(y)|^p\bigr|\,V_+(y)\,dy \\
&\qquad\qquad\le p\biggl(\int_{\Om}\bigl||\vp_k(y)|-|\vp(y)|\bigr|^p\, V_+(y)\,dy\biggr)^{\f{1}{p}}  \\
&\qquad\qquad\qquad\qquad\times\biggl(\int_\Om\bigl(|\vp_k(y)|\vee|\vp(y)|\bigr)^p\,V_+(y)\,dy\biggr)^{\f{p-1}{p}} \\
&\qquad\qquad\ls\|\vp\|_{\rX_0^{s,p}(\Om)}\|\vp_k-\vp\|_{\rX_0^{s,p}(\Om)}\|V_+\|_{L^q(\Om)}\to 0
\,\,\,\,\text{ as $k\to\iy$.}
\end{split}\end{equation}
Also, we claim that
\begin{equation}\lim_{k\to\iy}\int_{\BR^n}|\vp_k(y)|^p\,V_-(y)\,dy=\int_{\BR^n}|\vp(y)|^p\,V_-(y)\,dy;
\end{equation}
indeed, we have that 
$$V_-(y)\le[V_-,V_+^{b,1}]_{(A_1,A_1)}V_+^{b,1}(y)\,\,\text{ $\aee$ $y\in\BR^n$,}$$
because $(V_-,V_+^{b,1})\in(A_1,A_1)$  by (3.3). Since $V_+^b\in L^q_{\rm loc}(\BR^n)$ by (3.3), as in (3.5) we obtain that
\begin{equation*}\begin{split}
&\biggl|\int_{\BR^n}|\vp_k(y)|^p\,V_-(y)\,dy-\int_{\BR^n}|\vp(y)|^p\,V_-(y)\,dy\biggr| \\
&\qquad\qquad\le\int_{\BR^n}\bigl||\vp_k(y)|^p-|\vp(y)|^p\bigr|\,V_-(y)\,dy \\
&\qquad\qquad\le\bigl(\|b\|_{L^\iy(\BR^n)}\vee 1\bigr)[V_-,V_+^{b,1}]_{(A_1,A_1)} \\
&\qquad\qquad\qquad\qquad\times \int_{\Om}\bigl||\vp_k(y)|^p-|\vp(y)|^p\bigr|\,V^b_+(y)\,dy\to 0\,\,\,\, \text{ as $k\to\iy$.}
\end{split}\end{equation*}
Thus, by (3.2), (3.4) and (3.6), we obtain that
\begin{equation*}\begin{split}
\int_{\BR^n}|\vp(y)|^p\,V_-(y)\,dy&=\lim_{k\to\iy}\int_{\BR^n}|\vp_k(y)|^p\,V_-(y)\,dy \\
&\le\lim_{k\to\iy}\int_{\BR^n}|\vp_k(y)|^p\,V_+(y)\,dy=\int_{\BR^n}|\vp(y)|^p\,V_+(y)\,dy.
\end{split}\end{equation*}
Hence we are done. \qed

\begin{lemma} If $V\in\cP^{s,p}_q(\BR^n)$ for $q>\f{n}{ps}>1$ $($$\,p>1$, $0<s<1$$)$, then $\rY^{s,p}_0(\Om)$ is a quasi-Banach space with the quasinorm $\|\cdot\|_{\rY^{s,p}_0(\Om)}$ given by
$$\|u\|^p_{\rY^{s,p}_0(\Om)}:=\|u\|^p_{\rX^{s,p}_0(\Om)}+\int_\Om|u(y)|^p V(y)\,dy,\,\,u\in\rY^{s,p}_0(\Om).$$
Moreover, $\rY^{s,p}_0(\Om)=\rX^{s,p}_0(\Om)$ and they are quasinorm-equivalent.
\end{lemma}

\pf It follows from Lemma 3.1 and (2.5) that
\begin{equation*}\begin{split}\|u\|^p_{\rX_0^{s,p}(\Om)}&\le\|u\|^p_{\rY_0^{s,p}(\Om)}\le\|u\|^p_{\rX_0^{s,p}(\Om)}+\|u\|^p_{L^p_{V_+}(\Om)}  \\
&\le(1+|\Om|^{\f{ps}{n}-\f{1}{q}})\|u\|^p_{\rX_0^{s,p}(\Om)}.
\end{split}\end{equation*}
Since $\rX_0^{s,p}(\Om)$ is a Banach space, this implies the required results. \qed

\begin{thm} If $V\in\cP^{s,p}_q(\BR^n)$ for $q>\f{n}{ps}>1$ and $0<s<1$, then 
\begin{equation}\int_{\BR^n}|\vp(y)|\,V_-(y)\,dy\le\int_{\BR^n}|\vp(y)|\,V_+(y)\,dy
\end{equation}
for any $\vp\in\rY^{s,p}_0(\Om)$.
\end{thm}

\pf Take any $\vp\in\rY^{s,p}_0(\Om)$. Since $C^\iy_c(\Om)$ is dense in $\rX^{s,p}_0(\Om)$ 
(see \cite{G} and \cite{FSV}), $\rY^{s,p}_0(\Om)=\rX^{s,p}_0(\Om)$ and they are quasinorm-equivalent by Lemma 3.2, we can take a sequence $\{\vp_i\}\subset C^\iy_c(\Om)$ such that $\vp_i\to\vp$ in $\rY^{s,p}_0(\Om)$. So by (2.5) we also have that $\vp_i\to\vp$ in $L^{p_*}(\Om)$, where $p_*=\f{pn}{n-ps}$. So we can choose a subsequence $\{\vp_{i_k}\}$ such that 
\begin{equation}\vp_{i_k}\to\vp\,\,\text{ $\aee$ in $\Om$. }
\end{equation}
Also we have that 
\begin{equation}\lim_{k\to\iy}\int_{\BR^n}|\vp_{i_k}(y)|\,V_+(y)\,dy=\int_{\BR^n}|\vp(y)|\,V_+(y)\,dy;
\end{equation}
indeed, it follows from the fractional Sobolev inequality and H\"older's inequality that
\begin{equation*}\begin{split}
&\biggl|\int_{\BR^n}|\vp_{i_k}(y)|\,V_+(y)\,dy-\int_{\BR^n}|\vp(y)|\,V_+(y)\,dy\biggr| \\
&\quad\qquad\le\int_{\Om}\bigl||\vp_{i_k}(y)|-|\vp(y)|\bigr|\,V_+(y)\,dy \\
&\quad\qquad\le\biggl(\int_\Om\bigl||\vp_{i_k}(y)|-|\vp(y)|\bigr|^p\,V_+(y)\,dy\biggr)^{\f{1}{p}}\biggl(\int_\Om V_+(y)\,dy\biggr)^{\f{p-1}{p}} \\
&\quad\qquad\ls|\Om|^{\f{s}{n}+\f{1}{q'}-\f{1}{p}}\,\|\vp_{i_k}-\vp\|_{\rY^{s,p}_0(\Om)}\,\|V_+\|_{L^q(\Om)}
\to 0\,\,\text{ as $k\to\iy$ }
\end{split}\end{equation*}
where $q'$ is the dual exponent of $q$. 
Hence, by Fatou's lemma, Lemma 3.1, Lemma 3.5 in \cite{K2}, (3.2), (3.4), (3.8) and (3.9), we conclude that
\begin{equation*}\begin{split}
\int_{\BR^n}|\vp(y)|\,V_-(y)\,dy&=\lim_{p\to 1^+}\int_{\BR^n}|\vp(y)|^p\,V_-(y)\,dy \\
&\le\lim_{p\to 1^+}\liminf_{k\to\iy}\int_{\BR^n}|\vp_{i_k}(y)|^p\,V_-(y)\,dy \\
&\le\lim_{p\to 1^+}\liminf_{k\to\iy}\int_{\BR^n}|\vp_{i_k}(y)|^p\,V_+(y)\,dy \\
&\le\lim_{p\to 1^+}\int_{\BR^n}|\vp(y)|^p\,V_+(y)\,dy \\
&=\int_{\BR^n}|\vp(y)|\,V_+(y)\,dy.
\end{split}\end{equation*}
Therefore we complete the proof.
 \qed

\section{Local properties of weak subsolutions}

In this section, we shall obtain certain local properties for weak subsolutions of the
nonlocal $p$-Laplacian type Schr\"odinger equation. These results play a crucial role
in establishing a nonlocal H\"older continuity for weak solutions of the nonlocal equation (1.3).
In order to establish the result, we need several steps.

\begin{lemma} For $p>1$ and $M>0$, let $h(t)=t^{p-1}(t-M)_+-(t-M)_+^p\ge 0$. Then $h$ is Lipschitz continuous on $\BR$.
\end{lemma}

\pf We note that $h(t)=0$ for $t<M$, $h$ is in $C^1(M,\iy)$,
$$\lim_{t\to M^+}\f{h(t)-h(M)}{t-M}=0\,\,\text{ and }\,\,\lim_{t\to M^-}\f{h(t)-h(M)}{t-M}=M^{p-1}.$$
Also it is easy to check that
$\lim_{t\to\iy}h(t)=0$, because
$$\lim_{t\to\iy}\f{t^{p-1}(t-M)_+}{(t-M)_+^p}=1.$$
Moreover, in order to check the differentiability of $h$ at the infinity, set $g(t)=h(1/t)$.
Then we have that
$$\lim_{t\to 0^+}\f{g(t)-g(0)}{t}=\lim_{t\to 0^+}\biggl[\f{1}{t^p}\bigl(\f{1}{t}-M\bigr)-\f{1}{t}\bigl(\f{1}{t}-M\bigr)^p\biggr]=0,$$
because $$\lim_{t\to 0^+}\biggl[\f{1}{t^p}\biggl(\f{1}{t}-M\bigr)\bigl/\f{1}{t}\bigl(\f{1}{t}-M\bigr)^p\biggr]=\lim_{t\to 0^+}\f{1-tM}{(1-tM)^p}=1.$$
Thus this implies the required result. \qed

\begin{cor} If $M>0$ and $u\in\rX^{s,p}(\Om)$ for $p>1$ and $0<s<1$, then 

$(a)$ $|u|^{p-2}u(u-M)_+-(u-M)_+^p\in\rX^{s,p}(\Om)$ and it is nonnegative in $\BR^n$, and

$(b)$ $\bigl[|u|^{p-2}u(u-M)_+-(u-M)_+^p\bigr]\zt^p\in\rX_0^{s,p}(\Om)$ for any nonnegative $\zt\in C^\iy_c(\Om)$.
\end{cor}

\pf (a) Note that $|u|^{p-2}u(u-M)_+-(u-M)_+^p=h\circ u$ for the function $h$ given in Lemma 4.1. 
Since $h$ is Lipschitz continuous on $\BR$ by Lemma 4.1, it is obvious that $0\le h\circ u\in\rX^{s,p}(\Om)$.

(b) By the mean value theorem, we have that
\begin{equation*}\begin{split}
\bigl|\zt^p(x)-\zt^p(y)\bigr|&=\biggl|\int_{\zt(y)}^{\zt(x)}\f{d}{d\tau}\tau^p\,d\tau\biggr|
\le p\bigl(\zt^{p-1}(x)\vee\zt^{p-1}(y)\bigr)|\zt(x)-\zt(y)| \\
&\le 2p\,\|\zt\|^{p-1}_{L^\iy(\BR^n)}|\zt(x)-\zt(y)|
\end{split}\end{equation*} for all $x,y\in\BR^n$.
This means that $\zt^p$ is Lipschitz continuous on $\BR^n$. So we can easily derive the result.
Hence we complete the proof. \qed

\,

\,\,Next we need the following {\it nonlocal Caccioppoli type inequality} which is Theorem 1.5. This is a very useful tool in proving an interior H\"older regularity of weak solutions to the nonlocal $p$-Laplacian type Schr\"odinger equation.

\,

\,\,\,{\bf [Proof of Theorem 1.5.]} For simplicity of the proof, we assume that
$x_0=0$. 
Let $w=(u-M)_+$ for $M\in[0,\iy)$ and take any nonnegative $\zt\in
C_c^{\iy}(B_r)$. We use $\vp=w\zt^p$ as a testing function in
the weak formulation of the equation. Then we have that
\begin{equation}\begin{split}
\la u,\vp\ra_{H_p,K}+\int_{\BR^n}|u(y)|^{p-2}u(y)\vp(y)V(y)\,dy\le 0,
\end{split}\end{equation}
where $d_K(x,y)=K(x-y)\,dx\,dy$ and
$$\la u,\vp\ra_{H_p,K}=\iint_{\BR^n\times\BR^n}H_p(u(x)-u(y))(\vp(x)-\vp(y))\,d_K(x,y)\,\text{ for $p\ge 1$.}$$
The first term in the left-hand side of the above inequality can be decomposed into two parts as follows;
\begin{equation}\begin{split}\la u,\vp\ra_{H_p,K}&=\iint_{B_r\times B_r}H_p(u(x)-u(y))(\vp(x)-\vp(y))\,d_K(x,y) \\
&\quad\qquad+2\iint_{(\BR^n\s B_r)\times B_r}H_p(u(x)-u(y))\,\vp(x)\,d_K(x,y) \\
&:=I_1+2 I_2.
\end{split}\end{equation}
For the estimate of $I_1$, without loss of generality we may assume that $u(x)\ge u(y)$. Then we first observe that $w(x)\ge w(y)$ and
\begin{equation}H_p(u(x)-u(y))(\vp(x)-\vp(y))\ge(w(x)-w(y))^{p-1}(\vp(x)-\vp(y))
\end{equation}
whenever $x,y\in B_r$; indeed, it can easily be checked by considering three possible occasions (i) $u(x), u(y)>M$, (ii) $u(x)>M$, $u(y)\le M$, and (iii) $u(y)\le u(x)\le M$.
Also we observe that
\begin{equation}\begin{split}
|\zt(x)w(x)-\zt(y)w(y)|^p&\le 2^{p-1}|w(x)-w(y)|^p(\zt^p(x) +\zt^p(y)) \\
&\quad+2^{p-1}(w^p(x)+w^p(y))|\zt(y)-\zt(x)|^p.
\end{split}\end{equation}
By (b) of Lemma 2.4, (4.3) and (4.4), we have that
\begin{equation}\begin{split}
&H_p(u(x)-u(y))(\vp(x)-\vp(y)) \\
&\quad\qquad\ge\f{1}{4}|w(x)-w(y)|^p(\zt^p(x)+\zt^p(y)) \\
&\qquad\qquad-d_p(w^p(x)+w^p(y))|\zt(y)-\zt(x)|^p \\
&\qquad\quad\ge\f{1}{2^{p+1}}|\zt(x)w(x)-\zt(y)w(y)|^p \\
&\qquad\qquad-\bigl(\f{1}{4}+d_p\bigr)(w^p(x)+w^p(y))|\zt(y)-\zt(x)|^p.
\end{split}\end{equation}
Thus it follows from (4.5) that
\begin{equation}\begin{split}
I_1&\ge\f{1}{2^{p+1}}\iint_{B_r\times B_r}|\zt(x)w(x)-\zt(y)w(y)|^p\,d_K(x,y) \\
&\quad-2^p\bigl(\f{1}{4}+d_p\bigr)\iint_{B_r\times B_r}(w^p(x)+w^p(y))|\zt(y)-\zt(x)|^p\,d_K(x,y).
\end{split}\end{equation}
For the estimate of $I_2$, we note that
\begin{equation*}\begin{split}
H_p(u(x)-u(y))\vp(x)&\ge-(u(y)-u(x))^{p-1}_+(u(x)-M)_+\zt^p(x)\\
&\ge-(u(y)-M)^{p-1}_+(u(x)-M)_+\zt^p(x)\\
&=-w^{p-1}(y)w(x)\zt^p(x)
\end{split}\end{equation*}
and thus we have that
\begin{equation}\begin{split}
I_2&\ge-\iint_{(\BR^n\s B_r)\times B_r}w^{p-1}(y)w(x)\zt^p(x)\,d_K(x,y) \\
&\ge-\biggl(\,\sup_{x\in\supp(\zt)}\int_{\BR^n\s B_r}w^{p-1}(y)\,K(x-y)\,dy\biggr)\int_{B_r}w(x)\zt^p(x)\,dx.
\end{split}\end{equation}
Finally, we claim that
\begin{equation}\int_{\BR^n}|u(y)|^{p-2}u(y)\vp(y)V(y)\,dy\ge\int_{\BR^n}w^p(y)\zt^p(y)V(y)\,dy.
\end{equation}
This is equivalent to the following inequality
\begin{equation*}\begin{split}
&\int_{\BR^n}\bigl[|u(y)|^{p-2}u(y)w(y)-w^p(y)\bigr]\zt^p(y)V_+(y)\,dy \\
&\qquad\qquad\ge\int_{\BR^n}\bigl[|u(y)|^{p-2}u(y)w(y)-w^p(y)\bigr]\zt^p(y)V_-(y)\,dy,
\end{split}\end{equation*}
whose proof is just a direct application of Lemma 3.2, Theorem 3.3 and Corollary 4.2.
Hence the required inequality can be obtained from (4.1), (4.2), (4.6), (4.7) and (4.8). \qed

\,\,Next, we shall obtain the local boundedness of such weak subsolutions which is Theorem 1.1 and a relation between the nonlocal tail terms of the positive part and the negative part of weak solutions of the nonlocal $p$-Laplacian type Schr\"odinger equation (1.3) in the following theorems. 

\,

\,\,\,{\bf [Proof of Theorem 1.1.]} Take any $\zt\in C_c^{\iy}(B^0_r)$ such that $|\n\zt|\le c/r$ on $\BR^n$.
Let $w=(u-M)_+$ for $M\in(0,\iy)$.
By Lemma 3.2, Theorem 1.6 and the mean value theorem, we have that
\begin{equation}\begin{split}
&\iint_{B^0_r\times B^0_r}|\zt(x)w(x)-\zt(y)w(y)|^p\,d_K(x,y)  \\
&\qquad\qquad\qquad\ls\,r^{p-ps}\|\n\zt\|^p_{L^{\iy}(B^0_r)}\|w\|^p_{L^p(B^0_r)}+\cA(w,\zt,r,s)\,\|w\|_{L^1(B^0_r)}
\end{split}\end{equation} where $\cA(w,\zt,r,s)$ is the value given by
$$\cA(w,\zt,r,s)=\sup_{x\in \supp(\zt)}\int_{\BR^n\s B^0_r}w^{p-1}(y)\,K(x-y)\,dy.$$
Applying Proposition 2.1 to (4.9), we obtain that
\begin{equation}\begin{split}
\biggl(\,\fint_{B^0_r}|w\zt|^{p\gm}\,dx\biggr)^{\f{1}{\gm}}
&\,\,\ls\bigl(r^{p-ps}\|\n\zt\|^p_{L^{\iy}(B^0_r)}+r^{-ps}\bigr)r^{ps}\fint_{B^0_r}|w|^p \,dx \\
&\qquad\qquad+\cA(w,\zt,r,s)\, r^{ps}\fint_{B^0_r}w\,dx
\end{split}\end{equation} where $\gm=\ds\f{n}{n-ps}>1$.
For $k=0,1,2,\cdots$, we set
\begin{equation*}\begin{split}&r_k=(1+2^{-k})r,\,\,\,\,r^*_k=\f{r_k+r_{k+1}}{2},\\
&M_k=M+(1-2^{-k})M_*,\,\,\,\,M^*_k=\f{M_k+M_{k+1}}{2},
\end{split}\end{equation*}
$w_k=(u-M_k)_+$ and $w^*_k=(u-M^*_k)_+$ for a constant $M_*>0$ to be determined later.
In (4.9), for $k=0,1,\cdots,$ we choose a function $\zt_k\in C_c^{\iy}(B^0_{r^*_k})$ with $\zt_k|_{B^0_{r_{k+1}}}\equiv 1$ such that $0\le\zt_k\le 1$ and $$|\n\zt_k|\le c\,2^{k+2}/r\,\,\text{ in $\BR^n$.}$$
For $k=0,1,2,\cdots$, we set
$$N_k=\biggl(\,\fint_{B^0_{r_k}}|w_k|^p\,dx\biggr)^{\f{1}{p}}.$$
Since $w_k^*\ge w_{k+1}$ and $$w_k^*(x)\ge M_{k+1}-M_k^*=2^{-k-2}M_*$$ whenever $u(x)\ge M_{k+1}$, we then have that
\begin{equation}\begin{split}
N_{k+1}&\ls \biggl(\f{1}{|B^0_{r_k}|}\int_{B^0_{r_{k+1}}}\f{w^p_{k+1}(w_k^*)^{p(\gm-1)}}{(M_{k+1}-M_k^*)^{p(\gm-1)}}\,dx\biggr)^{\f{1}{p}} \\
&\ls \biggl(\f{2^k}{M_*}\biggr)^{\gm-1}\biggl(\,\fint_{B^0_{r_k}}|w^*_k\,\zt_k|^{p\gm}\,dx\biggr)^{\f{1}{p}}.
\end{split}\end{equation}
Since $\zt_k\in C_c^{\iy}(B^0_{r_k})$, $w_k^*\le w_0$ for all $k$ and
\begin{equation*}\begin{split}|y-x|&\ge|y-x_0|-|x-x_0| 
\ge\biggl(1-\f{r_k^*}{r_k}\biggr)|y-x_0|\ge 2^{-k-2}|y-x_0|
\end{split}\end{equation*}
for any $x\in B^0_{r^*_k}$ and $y\in\BR^n\s B^0_{r_k}$, we easily obtain that
\begin{equation}
\cA(w_k^*,\zt_k,r_k,s)\le c\,2^{k(n+ps)}\,r^{-ps}\,[\cT_r(w_0;x_0)]^{p-1}.
\end{equation}
Since $0\le w_k^*\le w_k$ and $w_k(x)\ge M_k^*-M_k=2^{-k-2}M_*$ if $u(x)\ge M_k^*$, it follows from (4.10), (4.11) and (4.12) that
\begin{equation*}\begin{split}
\biggl(\f{2^k}{M_*}\biggr)^{-\f{p(\gm-1)}{\gm}}N^{\f{p}{\gm}}_{k+1}
&\le c\, 2^{pk}\fint_{B^0_{r_k}}|w^*_k|^p\,dx \\
&\qquad+c\,2^{k(n+ps)}\biggl(\f{r_k}{r}\biggr)^{ps}\,[\cT_r(w_0;x_0)]^{p-1}\,\fint_{ B^0_{r_k}}w^*_k\,dx \\
&\le c\, 2^{pk}N_k^p+c[\cT_r(w_0;x_0)]^{p-1}\,2^{k(n+ps)}\fint_{B^0_{r_k}}\f{w^*_k w_k^{p-1}}{(M^*_k-M_k)^{p-1}}\,dx \\
&\le c\,\biggl(2^{pk}+ 2^{k(n+ps)}\biggl(\f{2^k}{M_*}\biggr)^{p-1}\,[\cT_r(w_0;x_0)]^{p-1}\biggr)\,N_k^p.
\end{split}\end{equation*}
Taking $M^*$ in the above so that
$$M_*\ge \dt\,\cT_r(w_0;x_0)\,\,\text{ for $\dt\in(0,1]$,}$$
we obtain that
\begin{equation*}\f{N_{k+1}}{M_*}\le d_0\,a^k\,\biggl(\f{N_k}{M_*}\biggr)^{1+\e}
\end{equation*}
where $d_0=c^{\f{\gm}{p}}\dt^{-\f{p-1}{p}\gm}>0$, $a=2^{\f{\gm}{p}(n+2s+p-1)+\f{ps}{n-ps}}>1$ and $\e=\gm-1>0$. 

If $N_0\le d_0^{-\f{1}{\e}}a^{-\f{1}{\e^2}}M_*,$ then we set
$$M_*=\dt\,\cT_r(w_0;x_0)+c_0\,\dt^{-\f{(p-1)n}{s p^2}}a^{\f{(n-ps)^2}{p^2 s^2}}N_0$$ where $c_0=c^{\f{n}{s p^2}}$.
By Lemma 2.3, we conclude that
\begin{equation*}\begin{split}
\sup_{B^0_r}u&\le M+M_*\\
&\le M+\dt\,\cT_r(w_0;x_0)+c_0\dt^{-\f{(p-1)n}{sp^2}}a^{\f{(n-ps)^2}{p^2 s^2}}\biggl(\,\fint_{B^0_{2r}}(u-M)_+^p\,dx\biggr)^{\f{1}{p}}.
\end{split}\end{equation*}
Hence, taking $M\downarrow 0$ in the above estimate, we obtain the required result. \qed

\,

\,\,\,The third one is a lemma which furnishes a relation between the nonlocal tails of the positive and negative part of weak solutions to the nonlocal $p$-Laplacian type Schr\"odinger equation.

\begin{lemma} Let $V\in\cP^{s,p}_q(\BR^n)$, $g\in W^{s,p}(\BR^n)$ for $q>\f{n}{ps}>1$ $($$\,p>1$, $s\in(0,1)$$)$.
If $u\in \rY^{s,p}_g(\Om)$ is a weak solution of the nonlocal $p$-Laplacian type Schr\"odinger equation $(1.3)$ such that $u\ge 0$ in $B^0_R\subset\Om$, then we have the estimate
\begin{equation*}\cT_r(u_+;x_0)\lesssim\bigl(1+\|V_+\|_{L^q(\Om)}\bigr)\,\sup_{B^0_r}u+\biggl(\f{r}{R}\biggr)^{\f{ps}{p-1}}\cT_R(u_-;x_0)
\end{equation*}
for any $r\in(0,R)$.
\end{lemma}

\pf  Without loss of generality, we may assume that $x_0=0$.
Let $M=\sup_{B_r}u$ and $\vp(x)=w(x)\zt^p(x)$ where $w(x)=u(x)-2M$ and $\zt\in C_c^{\iy}(B_{3r/4})$ is a function satisfying that $\zt|_{B_{r/2}}\equiv 1$, $0\le\zt\le 1$ and $|\n\zt|\le c/r$ in $\BR^n$. Then we have that
\begin{equation}\begin{split}
0&=\iint_{B_r\times B_r}H_p(u(x)-u(y))(\vp(x)-\vp(y))\,d_K(x,y) \\
&\quad+2\int_{\BR^n\s B_r}\int_{B_r}H_p(u(x)-u(y))(u(x)-2M)\zt^p(x)\,d_K(x,y) \\
&\quad+\int_{\BR^n}V(x)|u(x)|^{p-2}u(x)(u(x)-2M)\zt^p(x)\,dx \\
&:=J_1+J_2+J_3.
\end{split}\end{equation}
Since the fact that 
\begin{equation}-2M\le w(x):=u(x)-2M\le-M\end{equation}
for any $x\in B_r$, by (a) of Lemma 2.4 we have that
$$H_p(w(x)-w(y))(w(x)\zt^p(x)-w(y)\zt^p(y))\ge-c_p 4^p M^p(\zt(x)-\zt(y))^p$$
for any $x,y\in B_r$,
it follows from simple calculation that
\begin{equation}\begin{split}
J_1&\ge-c_p 4^p M^p\iint_{B_r\times B_r}(\zt(x)-\zt(y))^p\,d_K(x,y) \\
&\gtrsim-M^p r^{-ps}|B_r|.
\end{split}\end{equation} The lower estimate on $J_2$ can be splitted as follows;
\begin{equation*}\begin{split}
J_2&\ge 4\int_{\BR^n\s B_r}\int_{B_r}M(u(y)-M)^{p-1}_+\,\zt^p(x)\,d_K(x,y) \\
&\qquad-4M\int_{E_M}\int_{B_r}(u(x)-u(y))^{p-1}_+\,\zt^p(x)\,d_K(x,y) \\
&:=J_{2,1}-J_{2,2},
\end{split}\end{equation*} where $E_M=\{y\in\BR^n\s B_r:u(y)<M\}$.
Since $(u(y)-M)_+\ge u_+(y)-M$, it follows from (c) of Lemma 2.4 that
\begin{equation*}(u(y)-M)^{p-1}_+\ge b_p u_+^{p-1}(y)-M^{p-1}
\end{equation*} where $b_p=\mathbbm{1}_{(1,2]}(p)+2^{-(p-1)}\mathbbm{1}_{(2,\iy)}(p)$.
Thus the lower estimate on $J_{2,1}$ can be obtained as
\begin{equation}\begin{split}
J_{2,1}&\ge d_2 M r^{-ps}|B_r|\,\bigl[\cT_r(u_+;0)\bigr]^{p-1}-d_3 M^p r^{-ps}|B_r|
\end{split}\end{equation}
with universal constants $d_2,d_3>0$. If $x\in B_r$ and $y\in E_M$, then we observe that
\begin{equation*}\begin{split}(u(x)-u(y))^{p-1}_+&\le a_p\bigl(|u(x)-M|^{p-1}+|M-u(y)|^{p-1}\bigr)\\
&\le a_p M^{p-1}+a_p(M+u_-(y)-u_+(y))^{p-1}\\
&\le a_p M^{p-1}+a_p(M+u_-(y))^{p-1}\\
&\le a_p(1+a_p)M^{p-1}+a_p^2\,[u_-(y)]^{p-1}
\end{split}\end{equation*} where $a_p=\mathbbm{1}_{(1,2]}(p)+2^{p-1}\mathbbm{1}_{(2,\iy)}(p)$,
because $u_+(y)< M+u_-(y)$ for any $y\in E_M$.
Since $u_-(y)=0$ for all $y\in B_R$, the upper estimate on $J_{2,2}$ can thus be achieved by
\begin{equation}\begin{split}
J_{2,2}&\le 4 a_p(1+a_p) M^p\int_{\BR^n\s B_r}\int_{B_r}\zt^p(x)\,d_K(x,y)\\
&+4 a_p M\int_{\BR^n\s B_R}\int_{B_r}[u_-(y)]^{p-1}\zt^2(x)\,d_K(x,y)\\
&\le d_4 M^p r^{-ps}|B_r|+d_5 M R^{-ps}|B_r|\,\bigl[\cT_R(u_-;0)\bigr]^{p-1}
\end{split}\end{equation}
with universal constants $d_4,d_5>0$. Thus, by (4.16) and (4.17), we have that
\begin{equation}\begin{split}
J_2&\ge -d M^p r^{-ps}|B_r|-d M R^{-ps}|B_r|\,\bigl[\cT_R(u_-;0)\bigr]^{p-1}\\
&\quad+e M r^{-ps} |B_r|\,\bigl[\cT_r(u_+;0)\bigr]^{p-1}
\end{split}\end{equation}
where $d,e>0$ are some universal constants depending only on $n,s,\ld$ and $\Ld$.

Finally, it follows from (4.14), H${\ddot {\rm o}}$lder's inequality and the fractional Sobolev inequality that
\begin{equation}\begin{split}J_3&\ge-2 M^p\int_{\BR^n}V_+(x)\zt^p(x)\,dx \\
&\ge-2M^p \|V_+\|_{L^q(\Om)}\biggl(\,\int_{\Om}\zt^{pq'}(x)\,dx\biggr)^{\f{1}{q'}} \\
&\ge-2M^p \|V_+\|_{L^q(\Om)}\biggl(\,\int_{\Om}\zt^{\f{pn}{n-ps}}(x)\,dx\biggr)^{\f{n-ps}{n}}|\Om|^{\f{1}{q'}-\f{n-ps}{n}} \\
&\ge-2M^p\|V_+\|_{L^q(\Om)}|\Om|^{\f{1}{q'}-\f{n-ps}{n}}\bigl(r^{-ps}\|\zt\|^p_{L^p(B_r)}+[\zt]^p_{W^{s,p}(B_r)}\bigr) \\
&\gs-\|V_+\|_{L^q(\Om)} M^p r^{-ps}|B_r|
\end{split}\end{equation} where $q>\f{n}{ps}>1$ and $1<q'<\f{n}{n-ps}$ with 
$$\f{1}{q}+\f{1}{q'}=1.$$
Hence the estimates (4.13), (4.15), (4.18) and (4.19) give the required estimate. \qed

\,

\,\,\,Next we shall obtain the local boundedness for nonnegative weak subsolutions of the nonlocal equation (1.3) by employing Theorem 1.1 and Lemma 4.3. It is interesting that this estimate no longer depends on the nonlocal tail term, whose proof is pretty simple.

\begin{thm} Let $V\in\cP^{s,p}_q(\BR^n)$ and $g\in W^{s,p}(\BR^n)$ for $q>\f{n}{ps}>1$ $($\,$p>1$, $s\in(0,1)$\,$)$.
If $u\in \rY^{s,p}_g(\Om)$ is a nonnegative weak solution of the nonlocal $p$-Laplace type Schr\"odinger equation $(1.3)$, then we have the estimate
$$\sup_{B^0_r} u\le C\,\biggl(\,\fint_{B^0_{2r}}u^p(x)\,dx\biggr)^{\f{1}{p}}$$
for any $r>0$ with $B^0_{2r}\subset\Om$.
\end{thm}

\pf We choose some $\dt\in(0,1]$ so that $1-\dt d_0>0$ and take any $r>0$ with $B^0_{2r}\subset\Om$ where $$d_0=c_0(1+\|V_+\|_{L^q(\Om)})>0$$ for the universal constant $c_0>0$ given in Lemma 4.3. Then it follows from Theorem 1.1 and Lemma 4.3 that
\begin{equation*}\begin{split}
\sup_{B^0_r}u&\le\dt\,d_0\biggl[\,\sup_{B^0_r}u+\cT_{2r}(u^-;x_0))\biggr]+C_0\,\dt^{-\f{(p-1)n}{sp^2}}\biggl(\,\fint_{B^0_{2r}}u^p(x)\,dx\biggr)^{\f{1}{p}}
\end{split}\end{equation*}
Since $\cT_{2r}(u^-;x_0)=0$, we can easily derive the required result by taking
$$C=\f{C_0\,\dt^{-\f{(p-1)n}{sp^2}}}{1-\dt d_0}.$$
Hence we complete the proof. \qed

\,\,

\section{The logarithm of a weak solution is a locally bounded mean oscillation function }
In this section, we prove that the logarithm of a weak supersolution to the nonlocal equation (1.3) is a function with locally bounded mean oscillation. To do this, the following tool which is called the {\it fractional Poincar\'e inequality} is very useful.

Let $n\ge 1$, $p\ge 1$, $s\in(0,1)$ and $sp<n$. For a ball $B\subset\BR^n$, let $u_B$ denote the average of $u\in W^{s,p}(B)$ over $B$, i.e.
\begin{equation*}u_B=\fint_{\,\,\,B} u(y)\,dy.
\end{equation*}
Then it was shown in [BBM, MS] that
\begin{equation}\|u-u_B\|^p_{L^p(B)}\le\f{c_{n,p}(1-s)|B|^{\f{sp}{n}}}{(n-sp)^{p-1}}\,[u]^p_{W^{s,p}(B)}
\end{equation}
with a universal constant $c_{n,p}>0$ depending only on $n$ and $p$, which is usually very useful in getting the logarithmic estimate of weak supersolutions.
Of course, the logarithmic estimate could be obtained as in \cite{DKP1}, but we will not apply their approach to achieve it. Our method to realize the logarithmic estimate is easier and more simple than their method.

\,\,

\,\,\,{\bf [Proof of Theorem 1.2.]}  For simplicity, we set $x_0=0$. So, in what follows, we write $B_r:=B^0_r$ for $r>0$. Take any $r>0$ so that $B_{2r}\subset B_R$ where $B_R\subset\Om$. Consider a radial function $\zt\in C^{\iy}_c(B_{3r/2})$ with values in $[0,1]$ such that
$\zt|_{B_{r}}\equiv 1$, $\zt|_{\BR^n\s B_{2r}}\equiv 0$ and 
$$|\n\zt|\ls\f{1}{r}\,\,\text{ in $\BR^n$. }$$
We use the function $$\vp(x)=\f{\zt^p(x)}{u^{p-1}_{b}(x)}$$ as a testing function to the nonlocal $p$-Laplacian type Schr\"odinger equation (1.3), where $u_{b}(x)=u(x)+b$.
Then we have that
\begin{equation}\begin{split}0&\le\iint_{\BR^n\times\BR^n}H_p(u_b(x)-u_b(y))(\vp(x)-\vp(y))\,d_K(x,y)
+\int_{\BR^n}V(x) H_p(u(x))\vp(x)\,dx \\
&=\iint_{B_{2r}\times B_{2r}}H_p(u_b(x)-u_b(y))(\vp(x)-\vp(y))\,d_K(x,y) \\
&\quad\qquad+2\int_{\BR^n\s B_{2r}}\int_{B_{2r}}H_p(u_b(x)-u_b(y))\,\vp(x)\,d_K(x,y)  \\
&\quad\qquad+\int_{\BR^n}V(x)|u(x)|^{p-2}u(x)\vp(x)\,dx \\
&:=H(u,\vp)+I(u,\vp)+J(u,\vp).
\end{split}\end{equation}
Without loss of generality, we may assume that 
$$u_b(x)\ge u_b(y)$$ for the estimate $H(u,\vp)$; for, by symmetry, the other case $u_b(x)<u_b(y)$ can be treated in the exactly same way. 
Then we have two possible occasions;  (a) $u_b(x)\le 2 u_b(y)$, (b) $u_b(x)>2 u_b(y)$. 

\,

[Case (a) : $u_b(y)\le u_b(x)\le 2 u_b(y)$] By the mean value theorem, we note that
\begin{equation}\begin{split}
\zt(x)\ge\zt(y)\,\,&\Rightarrow\,\,\zt^p(x)-\zt^p(y)=p\int_{\zt(y)}^{\zt(x)}\tau^{p-1}d\tau\le p\zt^{p-1}(x)(\zt(x)-\zt(y)), \\
\zt(x)<\zt(y)\,\,&\Rightarrow\,\,\zt^p(x)-\zt^p(y)=p\int_{\zt(x)}^{\zt(y)}(-\tau^{p-1})d\tau\le p\zt^{p-1}(x)(\zt(x)-\zt(y)).
\end{split}\end{equation}
Then it follows from (5.3) that
\begin{equation}\begin{split}
\vp(x)-\vp(y)&=\f{\zt^p(x)-\zt^p(y)}{u_b^{p-1}(y)}+\zt^p(x)\,\biggl(\f{1}{u_b^{p-1}(x)}-\f{1}{u_b^{p-1}(y)}\biggr) \\
&\le\f{p\zt^{p-1}(x)(\zt(x)-\zt(y))}{u_b^{p-1}(y)} \\
&\quad+\zt^p(x)\int_0^1\f{d}{d\tau}\biggl(\f{1}{[\tau(u_b(x)-u_b(y))+u_b(y)]^{p-1}}\biggr)\,d\tau \\
&\le\f{p\zt^{p-1}(x)(\zt(x)-\zt(y))}{u_b^{p-1}(y)}  \\
&\quad\qquad-(p-1)\f{\zt^p(x)(u_b(x)-u_b(y))}{u_b^p(x)} \\
&\le\f{p\zt^{p-1}(x)\,|\zt(x)-\zt(y)|\,u_b(y)}{u_b^p(y)}  \\
&\quad\qquad-\f{(p-1)}{2^p}\,\f{\zt^p(x)(u_b(x)-u_b(y))}{u_b^p(y)}.
\end{split}\end{equation}
Applying Young's inequality with indices $p'=\f{p}{p-1}, p,\vep$, it follows from (5.2) that
\begin{equation}\begin{split}
&H(u,\vp)\le c_{n,p,s}\Ld p \\
&\iint_{B_{2r}\times B_{2r}}\f{\vep(u_b(x)-u_b(y))^p\zt^p(x)+c_\vep|\zt(x)-\zt(y)|^p u_b^p(y)}{u_b^p(y)} \f{dx\,dy}{|x-y|^{n+ps}}\\
&\qquad-\f{c_{n,p,s}\ld(p-1)}{2^p}\iint_{B_{2r}\times B_{2r}}\f{(u_b(x)-u_b(y))^p\zt^p(x)}{u_b^p(y)}\f{dx\,dy}{|x-y|^{n+ps}}.
\end{split}\end{equation}
If we choose $\vep=\f{\ld(p-1)}{2^{p+1}p\Ld}$ in (5.5), then we have that
\begin{equation}\begin{split}
H(u,\vp)&\ls-\iint_{B_{2r}\times B_{2r}}\zt^p(x)\f{(u_b(x)-u_b(y))^p}{u_b^p(y)}\f{\,dx\,dy}{|x-y|^{n+ps}} \\
&\quad\qquad+\iint_{B_{2r}\times B_{2r}}\f{|\zt(x)-\zt(y)|^p}{|x-y|^{n+ps}}\,dx\,dy \\
&\ls-\iint_{B_{2r}\times B_{2r}}\zt^p(x)\f{(u_b(x)-u_b(y))^p}{u_b^p(y)}\,\f{dx\,dy}{|x-y|^{n+ps}}+r^{n-ps}.
\end{split}\end{equation} because $x,y\in B_{2r}$.
Since $0\le u_b(x)-u_b(y)\le u_b(y)$, we have that
\begin{equation}\begin{split}
\bigl|\ln u_b(x)-\ln u_b(y)\bigr|^p&=\biggl(\int_0^1\f{u_b(x)-u_b(y)}{\tau(u_b(x)-u_b(y))+u_b(y)}\,d\tau\biggr)^p \\
&\le\f{(u_b(x)-u_b(y))^p}{u^p_b(y)}.
\end{split}\end{equation}
Thus by (5.6) and (5.7) we have that
\begin{equation}\begin{split}
H(u,\vp)&\ls-\iint_{B_{2r}\times B_{2r}}\zt^p(x)\,\biggl|\ln\,\biggl(\f{u_b(x)}{u_b(y)}\biggr)\biggr|^p\,\f{dx\,dy}{|x-y|^{n+ps}}+r^{n-ps} \\
&\ls-\iint_{B_{r}\times B_{r}}\,\biggl|\ln\,\biggl(\f{u_b(x)}{u_b(y)}\biggr)\biggr|^p\,\f{dx\,dy}{|x-y|^{n+ps}}+r^{n-ps}.
\end{split}\end{equation}

[Case (b) : $u_b(x)>2 u_b(y)$] It follows from the inequality in Lemma 2.5 with $\vep=(2^{p-1}-1)/2$ that
\begin{equation}\begin{split}
\vp(x)-\vp(y)&=\f{\zt^p(x)-\zt^p(y)}{u_b^{p-1}(x)}+\zt^p(y)\,\biggl(\f{1}{u_b^{p-1}(x)}-\f{1}{u_b^{p-1}(y)}\biggr) \\
&\le\f{\zt^p(x)-\zt^p(y)}{u_b^{p-1}(x)}+\zt^p(y)\,\biggl(\f{1}{2^{p-1}u_b^{p-1}(y)}-\f{1}{u_b^{p-1}(y)}\biggr) \\
&\le\f{\vep\zt^p(y)+c_\vep|\zt(x)-\zt(y)|^p}{u_b^{p-1}(x)}-(1-2^{-p+1})\f{\zt^p(y)}{u_b^{p-1}(y)} \\
&\le \f{c\, |\zt(x)-\zt(y)|^p}{u_b^{p-1}(x)}-\biggl(\f{1}{2}-\f{1}{2^p}\biggr)\,\f{\zt^p(y)}{u_b^{p-1}(y)}.
\end{split}\end{equation}
Since $u_b(x)\ge u_b(x)-u_b(y)\ge u_b(y),$ by (5.2) and (5.9) we have that
\begin{equation}\begin{split}
\f{H(u,\vp)}{c_{n,p,s}}&\le c\Ld\iint_{B_{2r}\times B_{2r}}\f{|\zt(x)-\zt(y)|^p}{|x-y|^{n+ps}}\,dx\,dy \\
&\quad-c\ld\biggl(\f{1}{2}-\f{1}{2^p}\biggr)\iint_{B_{2r}\times B_{2r}}\zt^p(y)\,\f{(u_b(x)-u_b(y))^{p-1}}{u_b^{p-1}(y)}\,\f{dx\,dy}{|x-y|^{n+ps}}.
\end{split}\end{equation} 
Since we know that $$(\ln t)^p\le c\,(t-1)^{p-1}\,\,\text{ for $t>2$,}$$ we have that
\begin{equation}\begin{split}
\bigl|\ln u_b(x)-\ln u_b(y)\bigr|^p&\le c\,\biggl(\f{u_b(x)-u_b(y)}{u_b(y)}\biggr)^{p-1} \\
&=c\,\f{(u_b(x)-u_b(y))^{p-1}}{u_b^{p-1}(y)}.
\end{split}\end{equation}
Combining (5.10) with (5.11), we have that
\begin{equation}\begin{split}
H(u,\vp)&\ls-\iint_{B_{2r}\times B_{2r}}\zt^p(y)\,\biggl|\ln\,\biggl(\f{u_b(x)}{u_b(y)}\biggr)\biggr|^p\,\f{dx\,dy}{|x-y|^{n+ps}} \\
&\qquad\qquad+\iint_{B_{2r}\times B_{2r}}\f{|\zt(x)-\zt(y)|^p}{|x-y|^{n+ps}}\,dx\,dy \\
&\ls-\iint_{B_{r}\times B_{r}}\,\biggl|\ln\,\biggl(\f{u_b(x)}{u_b(y)}\biggr)\biggr|^p\,\f{dx\,dy}{|x-y|^{n+ps}}+r^{n-ps}.
\end{split}\end{equation}
Hence, by (5.8) and (5.12), we conclude that
\begin{equation}H(u,\vp)\ls-\iint_{B_{r}\times B_{r}}\biggl|\ln\,\biggl(\f{u_b(x)}{u_b(y)}\biggr)\,\biggr|\,\f{dx\,dy}{|x-y|^{n+ps}}+r^{n-ps}.
\end{equation}

For the estimate of $I(u,\vp)$, we note that 
(i) $u(y)\ge 0$ and $u(x)-u(y)\le u(x)$ for $(x,y)\in B_{2r}\times(B_R\s B_{2r})$ and (ii) $(u(x)-u(y))_+\le u(x)+u_-(y)$ for $(x,y)\in B_{2r}\times(\BR^n\s B_R)$.
Since $\zt$ is supported in $B_{3r/2}$, the above observations (i) and (ii) yield that
\begin{equation}\begin{split}
I(u,\vp)&\le 2\,c_{n,p,s}\Ld\int_{B_{3r/2}}\int_{\BR^n\s B_{2r}}\f{1}{|y-x_0|^{n+ps}}\,dy\,dx \\
&\qquad+2\,c_{n,p,s}\Ld\int_{B_{3r/2}}\int_{\BR^n\s B_R}\f{u_-^{p-1}(y)}{u_b^{p-1}(x)}\,\f{dy\,dx}{|y-x_0|^{n+ps}} \\
&\ls r^{n-ps}+\f{r^{n-ps}}{b^{p-1}}\biggl(\f{r}{R}\biggr)^{ps} [\cT_R(u_-;x_0)]^{p-1}
\end{split}\end{equation}
because we see $|y-x|\ge|y-x_0|-|x-x_0|\ge|y-x_0|/4$ for any $(x,y)\in B_{3r/2}\times(\BR^n\s B_{2r})$.
Also, it follows from H\"older's inequality and Proposition 2.1 that
\begin{equation}\begin{split}
J(u,\vp)&\le\int_{\BR^n}V_+(x)\zt^p(x)\,dx \\
&\le\|V_+\|_{L^\tau(\Om)}\biggl(\,\int_\Om\zt^{p\tau'}(x)\,dx\biggr)^{\f{1}{\tau'}} \\
&\le\|V_+\|_{L^q(\Om)}\biggl(\,\int_{\Om}\zt^{\f{pn}{n-ps}}(x)\,dx\biggr)^{\f{n-ps}{n}}|\Om|^{\f{1}{q'}-\f{n-ps}{n}} \\
&\le\|V_+\|_{L^q(\Om)}|\Om|^{\f{1}{q'}-\f{n-ps}{n}}\bigl(r^{-ps}\|\zt\|^p_{L^p(B_r)}+[\zt]^p_{W^{s,p}(B_r)}\bigr) \\
&\le\|V_+\|_{L^q(\Om)}|\Om|^{\f{1}{q'}-\f{n-ps}{n}}r^{n-ps}\ls \|V_+\|_{L^q(\Om)}\,r^{n-ps}
\end{split}\end{equation}
where $q>\f{n}{ps}>1$, $p>1$ and $1<q'<\f{n}{n-ps}$ with 
$$\f{1}{q}+\f{1}{q'}=1.$$
By (5.13), (5.14) and (5.15), we obtain that
\begin{equation*}\begin{split}\iint_{B_r\times B_r}\,\biggl|\ln\biggl(\f{u(x)+b}{u(y)+b}\biggr)\biggr|^p&\,\f{dx\,dy}{|x-y|^{n+ps}}  \\
&\ls\f{r^{n-ps}}{b^{p-1}}\biggl(\f{r}{R}\biggr)^{ps}[\cT_R(u_-;x_0)]^{p-1}+r^{n-ps}\bigl(1+\|V_+\|_{L^q(\Om)}\bigr)
\end{split}\end{equation*}for any $b\in(0,1)$ and $r\in(0,R/2)$, since $x,y\in B_{2r}$.
Hence we complete the proof by applying (1.1). \qed





\,

\,\,\,We now introduce a kind of local BMO spaces on $B^0_R\subset\Om$, i.e. $\bmo^p(B^0_R)$ for $p>0$.
We define the norm $\|\cdot\|_{\bmo^p(B^0_R)}$ by
\begin{equation*}\|f\|_{\bmo^p(B^0_R)}=\sup_{r\in(0,R/2)}\biggl(\fint_{B^0_r}\bigl|f(y)-f_{B^0_r}\bigr|^p\,dy\biggr)^{\f{1}{p}}
\end{equation*} and the space 
$$\bmo^p(B^0_R)=\{f\in L^1_{\rm loc}(\BR^n): \|f\|_{\bmo^p(B^0_R)}<\iy\}.$$
When $p=1$, we write $\bmo^p(B^0_R)=\bmo(B^0_R)$.
Then we easily see that 
\begin{equation}\|\,|f|\,\|_{\bmo(B^0_R)}\le \|f\|_{\bmo(B^0_R)}
\end{equation} because $\bigl||f|-|f|_{B^0_r}\bigr|\le |f-f_{B^0_r}|$, and also
\begin{equation}\|f\pm g\|_{\bmo(B^0_R)}\le\|f\|_{\bmo(B^0_R)}+\|g\|_{\bmo(B^0_R)}.
\end{equation} We now observe that
$$a\wedge b=\f{a+b-|a-b|}{2}\,\,\text{ and }\,\,a\vee b=\f{a+b+|a-b|}{2}$$
for any $a,b\in\BR$. This implies that
\begin{equation}\begin{split}
\|f\vee g\|_{\bmo(B^0_R)}&\le\|f\|_{\bmo(B^0_R)}+\|g\|_{\bmo(B^0_R)},  \\
\|f\wedge g\|_{\bmo(B^0_R)}&\le\|f\|_{\bmo(B^0_R)}+\|g\|_{\bmo(B^0_R)}. 
\end{split}\end{equation}
In addition, we can obtain the following {\it John-Nirenberg inequality} (as in \cite{Gr}) by using the Calder\'on-Zygmund decomposition in harmonic analysis as follows; there exists some constants $b_1,b_2>0$ depending only on the dimension $n$ such that
$$\bigl|\{x\in B^0_r:|f(x)-f_{B^0_r}|>\ld\}\bigr|\le b_1\,e^{-(b_2/\|f\|_{\bmo(B_R^0)})\ld}|B^0_r|$$
for any $f\in\bmo(B^0_r)$, every $r>0$ with $B^0_{2r}\subset B^0_R$ and $B^0_R\subset\Om$, and every $\ld>0$. By standard analysis, this inequality makes it possible to easily show the following fact; 
\begin{equation}\begin{split}&\text{If $f\in\bmo(B^0_R)$ for $B^0_R\subset\Om$ and $1<p<\iy$, then $\|\cdot\|_{\bmo(B^0_R)}$ is } \\
&\,\,\text{ norm-equivalent to $\|\cdot\|_{\bmo^p(B^0_r)}$.}
\end{split}\end{equation}

\begin{lemma} If we set 
$$v(x)=\ds\ln\biggl(\f{a+b}{u(x)+b}\biggr)\,\,\text{ for $a,b\in(0,1)$ }$$ where the function $u$ satisfies the same assumption as Theorem 1.2, then we have that  
\begin{equation*}\fint_{B^0_r}|v(x)-v_{B^0_r}|^p\,dx\ls\,\fR_{b,r,R}(u_-;x_0)
\end{equation*} for any $r\in(0,R/2)$,
where $\fR_{b,r,R}(u_-;x_0)$ is given by
$$\fR_{b,r,R}(u_-;x_0)=\f{1}{b^{p-1}}\biggl(\f{r}{R}\biggr)^{ps}[\cT_R(u_-;x_0)]^{p-1}
+\bigl(1+\|V_+\|_{L^q(\Om)}\bigr).$$
It follows from this that $v\in\bmo(B^0_R)$ and moreover
$$\|v\|_{\bmo(B^0_R)}\le[\fR_{b,r,R}(u_-;x_0)]^{\f{1}{p}}<\iy.$$
\end{lemma}

\pf The first part easily follows from the fractional Poincar\'e inequality (5.1) and Theorem 1.2.
Also the second part can be shown by applying the Remark of Theorem 1.1 and H\"older's inequality because $u\in W^{s,p}(\BR^n)$. \qed

\begin{cor} If we set 
$$\bv=(v\vee 0)\wedge d\,\,\text{ for $d>0$ }$$ with the same $v$ as Lemma 5.1, then we have that 
\begin{equation*}\fint_{B^0_r}|\bv(x)-\bv_{B^0_r}|^p\,dx\ls\,\fR_{b,r,R}(u_-;x_0)
\end{equation*} for any $r\in(0,R/2)$, 
where $\fR_{b,r,R}(u_-;x_0)$ is given by
$$\fR_{b,r,R}(u_-;x_0)=\f{1}{b^{p-1}}\biggl(\f{r}{R}\biggr)^{ps}[\cT_R(u_-;x_0)]^{p-1}
+\bigl(1+\|V_+\|_{L^q(\Om)}\bigr).$$
It follows from this that $\bv\in\bmo(B^0_R)$ and moreover
$$\|\bv\|_{\bmo(B^0_R)}\le[\fR_{b,r,R}(u_-;x_0)]^{\f{1}{p}}<\iy.$$
\end{cor}

\pf  Without loss of generality, assume that $x_0=0$. By Lemma 5.1, we have that
\begin{equation*}\fint_{B_r}\bigl||v(x)|-|v|_{B_r}\bigr|^p\,dx\ls\,\fR_{b,r,R}(u_-;0),
\end{equation*}
because we see that
$$\bigl||v(x)|-|v|_{B_r}\bigr|\le |v(x)-v_{B_r}|.$$ Then we can easily derive from (5.17) and (5.19) that
\begin{equation*}
\fint_{B_r}|\bv(x)-\bv_{B_r}|\,dx\ls\,[\fR_{b,r,R}(u_-;0)]^{\f{1}{p}}.
\end{equation*}
Finally, the second part can be done as in Lemma 5.1. Hence we complete the proof. \qed

\section{Interior H$\ddot {\rm o}$lder regularity}

In this section, we establish an interior H\"older regularity of weak solutions to the nonlocal $p$-Laplacian type Schr\"odinger equation (1.3) by applying the previous results obtained in Section 4 and Section 5.

\,

\,\,\,{\bf [Proof of Theorem 1.3]} Fix any $p>1$ and $0<s<1$ and take any $R>0$ with $B_R(x_0)\subset\Om$.
For simplicity, without loss of generality, we may assume that $x_0=0$.
For any $k\in\BN\cup\{0\}$ and $r\in (0,R/2)$, we set 
\begin{equation*}r_k=\f{\dt^k r}{2}\text{ for $\dt\in\biggl(0,\bigl(\f{1}{4}\bigr)^{\f{p-1}{ps}}\biggr),$ $\,B_k=B_{r_k}$ and  $\, B^*_k=B_{2 r_k}$. }
\end{equation*}
Let us set 
\begin{equation*}\ds\Xi(r_0)=2\,\cT_{r/2}(u;0)+2\,C_0\,\biggl(\,\fint_{B_r}|u|^p\,dx\biggr)^{\f{1}{p}}
\end{equation*} where $C_0>1$ is the constant given in Theorem 1.1. As in the proof of Theorem 1.2, without loss of generality we may assume that 
$$\|u\|_{L^{\iy}(\BR^n)}=\f{1}{4C_0}<1.$$ So we see that 
\begin{equation*}0\le\Xi(r_0)<1.
\end{equation*}
For $k\in\BN\cup\{0\}$, we set
\begin{equation*}\Xi(r_k)=\biggl(\f{r_k}{r_0}\biggr)^{\e}\Xi(r_0)
\end{equation*}
where $\e\in (0,\f{ps}{p-1})$ is some constant to be determined later. For our proof, we have only to prove that
\begin{equation}\underset{B_k}\osc\,u\le\Xi(r_k)
\end{equation}
for any $k\in\BN\cup\{0\}$.

We proceed the proof by using the mathematical induction. By the remark of Theorem 1.1, we see that $$\underset{B_0}\osc\,u\le\Xi(r_0).$$ Assume that (6.1) holds for all $k\in\{0,1,\cdots, m\}$. Then we will show that (6.1) is still true for $m+1$. For this proof, we consider two possible cases; either
\begin{equation}\f{\bigl|B^*_{k+1}\cap\bigl\{u\ge\ds\inf_{B_k}u+\Xi(r_k)/2\bigr\}\bigr|}{|B^*_{k+1}|}\ge\f{1}{2}
\end{equation}
or
\begin{equation}\f{\bigl|B^*_{k+1}\cap\bigl\{u\le\ds\inf_{B_k}u+\Xi(r_k)/2\bigr\}\bigr|}{|B^*_{k+1}|}\ge\f{1}{2}\,.
\end{equation}
If (6.2) is true, then we set $u_k=u-\ds\inf_{B_k}u$, and if (6.3) is true, then we set $$u_k=\Xi(r_k)-\bigl(u-\inf_{B_k}u\bigr).$$ In these two cases, we see that $u_k\ge 0$ in $B_k$ and 
\begin{equation}\f{\bigl|B^*_{k+1}\cap\{u_k\ge\Xi(r_k)/2\}\bigr|}{|B^*_{k+1}|}\ge\f{1}{2}\,.
\end{equation} 
Furthermore, $u_k$ is a weak solution satisfying that
\begin{equation}\sup_{B_k}|u_m|\le 2\,\Xi(r_k)
\end{equation} for all $k\in\{0,1,\cdots,m\}$.
Under the induction hypothesis, if $m\ge 1$, then we now claim that
\begin{equation}[\cT_{r_k}(u_k;0)]^{p-1}\le c\,\dt^{-(p-1)\e}\,[\Xi(r_k)]^{p-1}
\end{equation} for all $k\in\{0,1,\cdots,m\}$.
Indeed, by (6.5) and the fact that
\begin{equation*}\begin{split}\int_{\BR^n\s B_0}\f{|u_m(x)|^{p-1}}{|x|^{n+ps}}\,dx &\ls r_0^{-ps}\sup_{B_0}|u|^{p-1} \\
&+r_0^{-ps}\,[\Xi(r_0)]^{p-1}+\int_{\BR^n\s B_0}\f{|u(x)|^{p-1}}{|x|^{n+ps}}\,dx \\
&\ls r_1^{-ps}\,\Xi(r_0),
\end{split}\end{equation*} 
we have the following estimate
\begin{equation*}\begin{split}
&[\cT_{r_m}(u_m;0)]^{p-1}=c\,r_m^{ps}\sum_{k=1}^m\int_{B_{k-1}\s B_k}\f{|u_m(x)|^{p-1}}{|x|^{n+ps}}\,dx +c\,r_m^{ps}\int_{\BR^n\s B_0}\f{|u_m(x)|^{p-1}}{|x|^{n+ps}}\,dx \\
&\qquad\quad\ls r_m^{ps}\sum_{k=1}^m\bigl[\,\sup_{B_{k-1}}|u_m|\,\bigr]^{p-1}\int_{\BR^n\s B_k}\f{1}{|x|^{n+ps}}\,dx +\,r_m^{ps}\int_{\BR^n\s B_0}\f{|u_m(x,t)|^{p-1}}{|x|^{n+ps}}\,dx \\
&\qquad\quad\le \sum_{k=1}^m\biggl(\f{r_m}{r_k}\biggr)^{ps}[\Xi(r_{k-1})]^{p-1}=[\Xi(r_0)]^{p-1}\sum_{k=1}^m\biggl(\f{r_m}{r_k}\biggr)^{ps}\biggl(\f{r_{k-1}}{r_0}\biggr)^{(p-1)\e} \\
&\qquad\quad=[\Xi(r_0)]^{p-1}\biggl(\f{r_m}{r_0}\biggr)^{(p-1)\e}\sum_{k=1}^m\biggl(\f{r_m}{r_k}\biggr)^{ps-(p-1)\e}\biggl(\f{r_{k-1}}{r_k}\biggr)^{(p-1)\e}       \\
&\qquad\quad=[\Xi(r_m)]^{p-1}\dt^{-(p-1)\e}\sum_{k=1}^{m}\dt^{(m-k)[ps-(p-1)\e]} \\
&\qquad\quad\le\dt^{-(p-1)\e}\,\f{\dt^{(m-1)[ps-(p-1)\e]}}{1-\dt^{-[ps-(p-1)\e]}}\,[\Xi(r_m)]^{p-1}\ls\dt^{-(p-1)\e}\,[\Xi(r_m)]^{p-1}. 
\end{split}\end{equation*}
For $k$ and $d>0$, we set 
\begin{equation*}v_k=\biggl[\ln\biggl(\f{\Xi(r_k)/2 +b}{u_k+b}\biggr)\vee 0\biggr]\wedge d.
\end{equation*}
Applying Corollary 5.2 with $a=\Xi(r_k)/2$, $b\in(0,1)$ and $d>0$, we have that
\begin{equation}\begin{split}\fint_{B^*_{k+1}}|v_k(x)-(v_k)_{B^*_{k+1}}|^p\,dx \ls\,\fR_{b,r_{k+1},r_k}((u_k)_-;0)
\end{split}\end{equation} where $\fR_{b,r,R}(u_-;x_0)$ is given by
$$\fR_{b,r,R}(u_-;x_0)=\f{1}{b^{p-1}}\biggl(\f{r}{R}\biggr)^{ps}[\cT_R(u_-;x_0)]^{p-1}
+\bigl(1+\|V_+\|_{L^q(\Om)}\bigr).$$
If we set $b=\dt^{\f{ps}{p-1}-\e}\,\Xi(r_k)$ in (6.7), then by (5.19) and (6.6) we obtain that
\begin{equation}\begin{split}&\fint_{B^*_{k+1}}|v_k(x)-(v_k)_{B^*_{k+1}}|\,dx  \\
&\quad\qquad\le\,\biggl(\,\fint_{B^*_{k+1}}|v_k(x)-(v_k)_{B^*_{k+1}}|^p\,dx\biggr)^{\f{1}{p}} \\
&\quad\qquad\le c\,\bigl(1+\|V_+\|_{L^q(\Om)}\bigr)^{\f{1}{p}}
\end{split}\end{equation}
where $c>0$ is the constant depending only on $n,s,p,\e,\ld$ and $\Ld$. 
From (6.4), we can derive the following estimate
\begin{equation}\begin{split}
d&=\f{1}{|B^*_{k+1}\cap\{u_k\ge\Xi(r_k)/2\}|}\int_{ B^*_{k+1}\cap\{u_k\ge\Xi(r_k)/2\}}d\,\,dx \\
&=\f{1}{|B^*_{k+1}\cap\{u_k\ge\Xi(r_k)/2\}|}\int_{B^*_{k+1}\cap\{v_k=0\}}d\,\,dx \\
&\le\f{2}{|B^*_{k+1}|}\int_{B^*_{k+1}}(d-v_k)\,dx= 2\bigl(d-(v_k)_{B^*_{k+1}}\bigr).
\end{split}\end{equation}
The estimates (6.8) and (6.9) make it possible to obtain the estimate
\begin{equation}\begin{split}
\f{|B^*_{k+1}\cap\{v_k=d\}|}{|B^*_{k+1}|}\,d&\le\f{2}{|B^*_{k+1}|}\int_{ B^*_{k+1}\cap\{v_k=d\}}\bigl(d-(v_k)_{B^*_{k+1}}\bigr)\,dx \\
&\le\f{2}{|B^*_{k+1}|}\int_{B^*_{k+1}\cap\{v_k=d\}}\bigl(v_k-(v_k)_{B^*_{k+1}}\bigr)\,dx \\
&\ls \bigl(1+\|V_+\|_{L^q(\Om)}\bigr)^{\f{1}{p}}. 
\end{split}\end{equation}
We now set 
\begin{equation*}d=d_*:=\ln\,\biggl(\f{\Xi(r_k)/2+\dt^{2s-\gm}\,\Xi(r_k)}{3\,\dt^{2s-\gm}\,\Xi(r_k)}\biggr).
\end{equation*}
Then we see that $d_*\sim\ln(1/\dt)$. By (6.10), we have that
\begin{equation}\begin{split}
\f{|B^*_{k+1}\cap\{u_k\le 2\,\dt^{\f{ps}{p-1}-\e}\,\Xi(r_k)\}|}{|B^*_{k+1}|}&\le\f{c\,\bigl(1+\|V_+\|_{L^q(\Om)}\bigr)^{\f{1}{p}}}{d_*} \\
&\le\f{c_0\,\bigl(1+\|V_+\|_{L^q(\Om)}\bigr)^{\f{1}{p}}}{\ln(1/\dt)}.
\end{split}\end{equation}

Now we proceed the next step with a well-known iteration process as follows.
For $i\in\BN\cup\{0\}$, we set
\begin{equation*}\rho_i=(1+2^{-i})r_{k+1},\, \bar\rho_i=\f{\rho_i+\rho_{i+1}}{2},\,B_i=B_{\rho_i}\,\text{ and }\,{\bar B}_i=B_{\bar\rho_i}.
\end{equation*}
For $i\in\BN\cup\{0\}$, we consider a function $\zt_i\in C^{\iy}_c(B_{\bar\rho_i})$ with $\zt_i|_{B_{\rho_{i+1}}}\equiv 1$ such that $0\le\zt_i\le 1$ and $|\n\zt_i|\le c\rho_i^{-1}$ in $\BR^n$. 
Moreover we set 
$d_i=(1+2^{-i})\dt^{\f{ps}{p-1}-\e}\,\Xi(r_k)$ and $w_i=(d_i-u_k)_+$ and 
\begin{equation*}N_i=\f{|B_i\cap\{u_k\le d_i\}|}{|B_i|}=\f{|B_i\cap\{w_i\ge 0\}|}{|B_i|}.
\end{equation*}
From (6.11), we see that 
\begin{equation}N_0\le \f{c_0\bigl(1+\|V_+\|_{L^q(\Om)}\bigr)^{\f{1}{p}}}{\ln(1/\dt)}.
\end{equation}
By Theorem 1.5, we have that
\begin{equation*}\begin{split}
[w_i \zt_i]^p_{W^{s,p}(B_{\rho_i})}&\ls 
\iint_{B_{\rho_i}\times B_{\rho_i}}[w_i(x)\vee w_i(y)]^p|\zt_i(x)-\zt_i(y)|^p\,d_K(x,y)\\
&\qquad+\biggl(\,\sup_{x\in\supp(\zt_i)}\int_{\BR^n\s B_{\rho_i}}[w_i(y)]^{p-1}\,K(x-y)\,dy\biggr)\|w_i\zt_i^p\|_{L^1(B_i)} \\
&\quad:=A(\rho_i,w_i,\zt_i)+B(\rho_i,w_i,\zt_i).
\end{split}\end{equation*}
Then we have the following estimate
\begin{equation}\begin{split}
A(\rho_i,w_i,\zt_i)&\ls\,d_i^p\int_{B_{\rho_i}}\int_{B_{\rho_i}\cap\{u_k\le d_i\}}
\f{\sup_{\BR^n}|\n\zt_i|^p}{|x-y|^{n+ps-p}}\,dx\,dy \\
&\ls\,d_i^p\biggl(\f{1}{\rho_i}\biggr)^p\int_{B_{\rho_i}\cap\{u_k\le d_i\}}\int_{B_{2\rho_i}}\f{1}{|y|^{n+2s-2}}\,dy\,dx \\
&\ls\,\,d_i^p\, \rho_i^{-ps}\,|B_i\cap\{u_k\le d_i\}|.
\end{split}\end{equation}
From the fact that
$$|y-x|\ge|y|-|x|\ge\biggl(1-\f{\bar\rho_i}{\rho_i}\biggr)|y|\ge 2^{-i-2}|y|$$
for any $y\in\BR^n\s B_{\rho_i}$ and $x\in B_{\bar\rho_i}$, we obtain that
\begin{equation}\begin{split}
B(\rho_i,w_i,\zt_i)&\ls\,d_i\,2^{i(n+ps)}\,|B_i\cap\{u_k\le d_i\}|\int_{\BR^n\s B_{\rho_i}}\f{|w_i(y)|^{p-1}}{|y|^{n+ps}}\,dy \\
&\ls\,2^{i(n+ps)}\,d_i\,\rho_i^{-ps}\,|B_i\cap\{u_k\le d_i\}|\,[\cT_{r_{k+1}}(w_i;0)]^{p-1}.
\end{split}\end{equation}
Thus it follows from (6.13) and (6.14) that
\begin{equation}\begin{split}
&[w_i\zt_i]^p_{W^{s,p}(B_{\rho_i})} \\
&\quad\ls\,\bigl(\,d_i^p+2^{i(n+ps)}\,d_i\,[\cT_{r_{k+1}}(w_i;0)]^{p-1}\bigr)\rho_i^{-ps}\,|B_i\cap\{u_k\le d_i\}|.
\end{split}\end{equation}
From (6.6) and the facts that $w_i\le 2\dt^{\f{ps}{p-1}-\e}\,\Xi(r_k)$ in $B_k$ and $w_i\le|u_k|+2\dt^{\f{ps}{p-1}-\e}\,\Xi(r_k)$ in $\BR^n$, we can derive that
\begin{equation}\begin{split}
&[\cT_{r_{k+1}}(w_i;0)]^{p-1} \\
&\qquad\ls\,r^{ps}_{k+1}\int_{B_k\s B_{k+1}}\f{|w_i(y,t)|^{p-1}}{|y|^{n+ps}}\,dy +\biggl(\f{r_{k+1}}{r_k}\biggr)^{ps}[\cT_{r_k}(w_i;0)]^{p-1} \\
&\qquad\ls\,\dt^{ps-(p-1)\e}\,[\Xi(r_k)]^{p-1}+\dt^{ps}[\cT_{r_k}(u_k;(0,0))]^{p-1}\ls\,d_i^{p-1}.
\end{split}\end{equation}
Thus by (6.15) and (6.16), we have that
\begin{equation}\begin{split}[w_i\zt_i]^p_{W^{s,p}(B_{\rho_i})}\ls\,2^{i(n+ps)}\,d^p_i\,\rho_i^{-ps}\,|B_i\cap\{u_k\le d_i\}|.
\end{split}\end{equation}
By applying (6.17) and the fractional Sobolev's inequality with the exponent
$$\gm=\f{n}{n-ps},$$
we can deduce the following inequalities
\begin{equation}\begin{split}
\biggl(\int_{B_{i+1}}|w_i|^{p\gm}\,dx\biggr)^{\f{1}{\gm}}&\le\biggl(\int_{B_i}|w_i\zt_i|^{\f{pn}{n-ps}}\,dx\biggr)^{\f{n-ps}{n}} \\
&\ls[w_i\zt_i]^p_{W^{s,p}(B_{\rho_i})}+\rho^{-ps}_i\|w_i\zt_i\|^p_{L^p(B_{\rho_i})} \\
&\ls\,2^{i(n+ps)}\,d^p_i\,\rho_i^{-ps}\,|B_i\cap\{u_k\le d_i\}|.
\end{split}\end{equation}
Since $|B_{i+1}|\sim\rho_i^{n}\sim|B_i|$ and
$$w_i=(d_i-u_k)_+\ge(d_i-d_{i+1})\,\mathbbm{1}_{\{u_k\le d_{i+1}\}}\ge 2^{-i-2}d_i\,\mathbbm{1}_{\{u_k\le d_{i+1}\}},$$
this estimate (6.18) yields that
\begin{equation*}\begin{split}
(d_i-d_{i+1})^p\biggl(\f{|B_{i+1}\cap\{u_k\le d_{i+1}|}{|B_{i+1}|}\biggr)^{\f{1}{\gm}}
&\le\f{\rho_i^{ps}}{|B_i|}\biggl(\int_{B_{i+1}}|w_i|^{p\gm}\,dx\,dt\biggr)^{\f{1}{\gm}} \\
&\ls\,2^{i(n+ps)}\,d^p_i\,\f{|B_i\cap\{u_k\le d_i\}|}{|B_{i}|},
\end{split}\end{equation*}
which gives that
$$N_{i+1}^{\f{1}{\gm}}\le c\,\f{2^{i(n+ps)}\,d_i^p}{(d_i-d_{i+1})^p}\,N_i\le c\,2^{i(n+ps+p)}N_i.$$
This leads us to obtain that
\begin{equation*}N_{i+1}\le c_1\,2^{i\gm(n+ps+p)}\,N_i^{1+\f{ps}{n-ps}}.
\end{equation*}
If we could show that
\begin{equation}N_0=\f{|B_0\cap\{u_k\le 2\,\dt^{\f{ps}{p-1}-\e}\,\Xi(r_k)\}|}{|B_0|}\le c_1^{-\f{n-ps}{ps}}2^{-\f{\gm (n-ps)^2(n+ps+p)}{p^2 s^2}}:=c_*,
\end{equation}
then by Lemma 2.2 we conclude that $\lim_{i\to\iy}N_i=0$, i.e. 
$$\inf_{B_{k+1}}u_k>\dt^{\f{ps}{p-1}-\e}\,\Xi(r_k).$$
In order to guarantee (6.19), by (6.12) we observe that
\begin{equation*}c_*\ge\f{c_0\bigl(1+\|V_+\|_{L^q(\Om)}\bigr)^{\f{1}{p}}}{\ln(1/\dt)}\,\,\,\Leftrightarrow\,\,\,\dt\le e^{-(c_0/c_*)(1+\|V_+\|_{L^q(\Om)})^{1/p}},
\end{equation*}
and so we choose $\dt>0$ as 
$$\dt=e^{-(c_0/c_*)(1+\|V_+\|_{L^q(\Om)})^{1/p}}\wedge\bigl(\f{1}{4}\bigr)^{\f{p-1}{ps}}.$$
If $u_k=u-\inf_{B_k} u$, then by (6.1) we have that
\begin{equation}\underset {B_{k+1}}\osc\, u=\underset{B_{k+1}}\osc\, u_k
\le\underset{B_k}\osc\, u-\inf_{B_{k+1}}u_k\le(1-\dt^{\f{ps}{p-1}-\e})\,\Xi(r_k).
\end{equation}
If $u_k=\Xi(r_k)-(u-\inf_{B_k}u)$, then we have that
\begin{equation}\underset{B_{k+1}}\osc\, u=\underset{B_{k+1}}\osc\, u_k
=\Xi(r_k)-\inf_{B_{k+1}}u+\inf_{B_k}u-\inf_{B_{k+1}}u_k\le(1-\dt^{\f{ps}{p-1}-\e})\,\Xi(r_k).
\end{equation}
From (6.20) and (6.21), we obtain that
\begin{equation}\begin{split}\underset{B_{k+1}}\osc\, u&\le(1-\dt^{\f{ps}{p-1}-\e})\,\Xi(r_k)
=(1-\dt^{\f{ps}{p-1}-\e})\biggl(\f{r_k}{r_{k+1}}\biggr)^{\e}\,\Xi(r_{k+1}) \\
&=(1-\dt^{\f{ps}{p-1}-\e})\dt^{-\e}\,\Xi(r_{k+1}).
\end{split}\end{equation}
To find $\e$ so that $(1-\dt^{\f{ps}{p-1}-\e})\dt^{-\e}\le 1$, we consider the function
$$\xi(\e)=\dt^\e+\dt^{\f{ps}{p-1}-\e}.$$
We note that 
$$\xi'(\e)=\dt^\e\ln\dt\,\bigl(1-\dt^{\f{ps}{p-1}-2\e}\bigr)=0\,\,\,\,\Leftrightarrow
\,\,\,\,\e=\f{ps}{2(p-1)}:=\e_{p,s}\,,$$
$\xi'(\e)<0$ for $\e<\e_{p,s}$ and $\xi'(\e)>0$ for $\e>\e_{p,s}$. These facts imply that the graph of $\xi$ is going down from the point $(0,1+\dt^{2\e_{p,s}})$ to the point $(\e_{p,s},1+2\dt^{\e_{p,s}})$ and is going up the point $(2\e_{p,s},1+\dt^{2\e_{p,s}})$ right after that. Since we see that 
$$2\dt^{\e_{p,s}}<1+\dt^{2\e_{p,s}}\,\,\,\,\text{ and }\,\,\,\,2\dt^{\e_{p,s}}<1,$$ 
we can find exactly two $\e$'s inside $(0,\f{ps}{p-1})$ so that
\begin{equation}\dt^\e+\dt^{\f{ps}{p-1}-\e}=1.
\end{equation}
If we set $Y=\dt^\e$, then the above equation (6.23) will be transformed into
$$\f{Y^2-Y+\dt^{\f{ps}{p-1}}}{Y}=0$$
and its solutions are
$$\dt^\e=Y=\f{1\pm\sqrt{1-4\dt^{\f{ps}{p-1}}}}{2}\in(0,1),$$
because $\dt<(1/4)^{\f{p-1}{ps}}$. Hence it turns out that the solutions of (6.23) are
$$\e_0^{\pm}=\f{\ln\biggl(\,\ds\f{1\pm\sqrt{1-4\dt^{\f{ps}{p-1}}}}{2}\,\biggr)}{\ln\dt}.$$
Then we see that $\xi(\e)\ge 1$, i.e. $(1-\dt^{\f{ps}{p-1}-\e})\dt^{-\e}\le 1$ for all $\e\in(0,\e_0^-]\cup[\e_0^+,2\e_{p,s})$. Thus, if $\e\in(0,\e_0^-]\cup[\e_0^+,2\e_{p,s})$, then by (6.22) we conclude that 
$$\underset{B_{k+1}}\osc\,u\le\Xi(r_{k+1}).$$
Therefore we complete the proof. \qed

\,

\noindent{\bf Acknowledgement.} 
Yong-Cheol Kim was supported by School of Education, Korea University Grant in 2023.


\end{document}